\documentclass[12pt]{article}
\newtheorem{proposition}{{\bf Proposition}}
\newtheorem{lemma}{{\bf Lemma}}
\newtheorem{theorem}{{\bf Theorem}}
\newtheorem{corollary}{{\bf Corollary}}

\usepackage[psamsfonts]{amssymb}
\usepackage{amsfonts}
\usepackage[english]{babel}
\usepackage{amsmath}
\usepackage{a4wide}
\usepackage{graphicx}
\newcommand{\Rset}{\ensuremath{\mathbb{R}}}
\def\div{\mbox{{\rm div}}}
\def\ln{\mbox{{\rm ln}}}

\def\div{\mbox{{\rm div}}}

\date{}

\begin{document}

%\frontmatter

%%%%%%%%%%%%%%%%%%%%%%%%%%%%%%%%%%%%%%%%%%%%%%%%%%%%%%%%%%%%%%%%%%%%%%%%%%

\title{On global solutions to some non-Markovian quantum kinetic models of Fokker-Planck type}

\author{Miguel A. Alejo,\thanks{Departamento de Matem\'atica, Federal University of Santa Catarina, 88040-900 Florian\'opolis, BRAZIL. miguel.alejo@ufsc.br} \ Jos\'e Luis L\'opez\thanks{Departamento de
Matem\'atica Aplicada and Excellence Research Unit "Modeling Nature" (MNat), Facultad de Ciencias, Universidad de Granada,
18071 Granada, SPAIN.  jllopez@ugr.es}}

\date{}

\maketitle

%%%%%%%%%%%%%%%%%%%%%%%%%%%%%%%%%%%%%%%%%%%%%%%%%%%%%%%%%%%%%%%%%%%%%%%%%%

\begin{abstract}
In this paper, global well-posedness of the non-Markovian Unruh-Zurek and Hu-Paz-Zhang master equations with nonlinear electrostatic coupling is demonstrated. They both consist of a Wigner-Poisson like equation subjected to a dissipative Fokker-Planck mechanism with time-dependent coefficients of integral type, which makes necessary to take into account the full history of the open quantum system under consideration to describe its present state. From a mathematical viewpoint this feature makes particularly elaborated the calculation of the propagators that take part of the corresponding mild formulations, as well as produces rather strong decays near the initial time ($t=0$) of the magnitudes involved, which would be reflected in the subsequent derivation of a priori estimates and a significant lack of Sobolev regularity when compared with their Markovian counterparts. The existence of local-in-time solutions is deduced from a Banach fixed point argument, while global solvability follows from appropriate kinetic energy estimates.
\end{abstract}

%%%%%%%%%%%%%%%%%%%%%%%%%%%%%%%%%%%%%%%%%%%%%%%%%%%%%%%%%%%%%%%%%%%

\vskip 12pt

\noindent
{\bf AMS Subject classification: } 

\vskip 12pt


%%%%%%%%%%%%%%%%%

\section{Introduction}

Open quantum systems are on the theoretical basis of many dissipative and/or diffusive physical processes and experiments \cite{BLPV,BP,Dav,VA,Fr,RH,We}. They roughly consist of a quantum system interacting with a thermal reservoir. In the Wigner function representation, this kind of systems are typically described by a Fokker-Planck kernel that incorporates damping and quantum diffusion terms to the standard Liouvillian evolution. Indeed, if one considers the simplest case in which a particle  interacts with an idealized thermal bath consisting of an infinite number of harmonic oscillators, then the resulting quantum Fokker-Planck equation (once the degrees of freedom of the thermal bath have been traced out) ruling the evolution of the Wigner function $w(t,x,\xi)$ associated with the reduced density matrix $\rho(t,x,y)$, namely 
\begin{equation}
w(t, x,\xi)= \frac{1}{(2\pi)^3} \int_{\Rset^3} \rho\Big(t,x + \frac{\eta}{2}, x - \frac{\eta} {2}\Big) \, e^{-i \eta \cdot \xi} \, d\eta \, , 
\label{w}
\end{equation}
reads as follows: 
\begin{equation}
\partial_t w + (\xi \cdot \nabla_x) w  - 2 \gamma {\rm div}_\xi (\xi w) = \frac{D_{pp}}{m^2} \Delta_\xi w + \frac{2 D_{pq}}{m} \div_x (\nabla_\xi w) + D_{qq} \Delta_x w \, ,
\label{wfp}
\end{equation}
where the coefficients $\gamma$, $D_{pp}$, $D_{pq}$, and $D_{qq}$ are all positive constants related to the interaction. Here, $t>0$ is the time variable, while $x, \xi \in \Rset^3$ hold for the position and momentum of the particle, respectively.  It is important to notice that in this approximation $w(t,x,\xi)$ evolves in a Markovian way, in the sense that it suffices to know the initial condition $w(0,x,\xi)$ so as to predict the state of the system at any future time $t>0$. This model was introduced by Caldeira and Leggett in \cite{CL} as a kinetic description of quantum Brownian motion (see also \cite{Dio1,Dio2}), and several versions of it were analyzed from a mathematical point of view in \cite{ALMS, CLN, LN}, among other works. In what follows we are concerned with a nonlinear correction of two variants of Eq. (\ref{wfp}), accounting for the effects of the three-dimensional electrostatic Poisson potential, namely 
$$
 V(t,x) = \frac{1}{4\pi|x|} \ast_x n(t,x) \, ,
 $$
where $\ast_x$ denotes the convolution product in the position variable and 
\begin{equation}
n(t,x) = \int_{\Rset^3} w(t,x,\xi) \, d\xi
\label{n}
\end{equation}
stands for the position density. Its presence in Eq. (\ref{wfp}) (as well as in its variants) clearly makes the equation nonlinear through the pseudo-differential operator
\begin{eqnarray}
\nonumber
\Theta[V]w(t,x,\xi) &=& \frac{i}{(2\pi)^3} \int_{\Rset^6} \left(V\Big(t,x + \frac{\hbar \eta}{2m}\Big) - V\Big(t, x -\frac{\hbar \eta}{2m} \Big)\right)
w(t,x,\xi') e^{-i(\xi-\xi') \cdot \eta} \, d(\xi', \eta) \\ &=& H \ast_\xi w \, ,
\label{pseudo}
\end{eqnarray}
with
\begin{equation}
H(t,x,\xi) = \frac{16 m^3}{\hbar^3} {\hbox{\rm{Re}}}\left\{  i e^{\frac{2m}{\hbar}i x \cdot \xi} {\check{V}}\Big(t,\frac{2m}{\hbar}\xi \Big) \right\} \, ,
\label{H}
\end{equation}
${\check{V}}$ denoting the inverse Fourier transform ${\check{V}}(t,\xi) = (2\pi)^{-3} \int_{\Rset^3} V(t,x) e^{ix\cdot\xi} \, dx$.

Our purpose in this paper is to extend the well-posedness results already known for Eq. (\ref{wfp}) subjected to the nonlinearity (\ref{pseudo}) (that is, existence and regularity of a unique global mild solution to the initial value problem in an appropriate function space) to the non-Markovian regime, in which case it is required to know the whole past history of the system in order to describe its current state. At the level of the governing equations, this is translated into the fact that the coefficients will not be constant any more, as in the Caldeira-Leggett model, but time dependent via integrals of the type 
\begin{equation}
\int_0^t \sigma(s) \, ds  \, . 
\label{nm}
\end{equation}
In this situation, the calculation of the associated propagators is rather more delicate than for Eq. (\ref{wfp}), since one has to make sure that eventual divergences of the integrals involved will not show up. Besides, we are faced with the case in which the coefficient $D_{qq}$ is set to zero, generating an obvious lack of elliptic regularity in the position variable that gives rise to much worse estimates than those available for Eq. (\ref{wfp}) (see \cite{LN}). On that basis, our systems of interest in this paper are the following two, usually known as Unruh-Zurek and Hu-Paz-Zhang master equations, respectively:
\begin{eqnarray*}
&&T[w] - \Omega_0^2 (x \cdot \nabla_\xi) w - 2 \gamma \div_{\xi}(\xi w) + \Theta[V]w = c(t) \Delta_{\xi} w - d(t) \div_x(\nabla_\xi w) \, , \\
&&T[w] - (\Omega^2 - 2 a(t)) (x \cdot \nabla_\xi) w - 2 b(t) \div_{\xi}(\xi w) + \Theta[V]w =  c(t) \Delta_{\xi} w - d(t) \div_x(\nabla_\xi w) \, ,
\end{eqnarray*}
where $T[w] = \partial_t w + (\xi \cdot \nabla_x) w$ is the transport operator and where the various (constant or time-dependent) coefficients are defined in the following section. Two main structural differences between them are noticeable. On one hand, the coefficient in front of the term $(x \cdot \nabla_\xi)w$ (namely, the oscillation frequency of the system) in the second equation is shifted by a time-dependent function of the type represented in (\ref{nm}), if compared to that of the first equation; on the other hand, the coefficient in front of the term $\div_{\xi}(\xi w)$ (that is, the dissipation rate of the system into the bath) is time-dependent (again of the type described in (\ref{nm})) in the second equation and constant in the first one.

Quantum non-Markovian environments have been shown to play an important role in the correct understanding of a wide variety of processes and systems stemming from quantum optics, quantum chemistry, biophysics or superconductor theory. In recent years, much attention has been paid to the description of memory effects in quantum channels that may increase their capacities and potentially contribute to the development of quantum information technologies, cryptography and teleportation \cite{MOP}. Furthermore, non-Markovian effects have proved fundamental to sustain quantum coherence of biomolecular excitons in photosynthetic complexes over long times, opening the way to the possibility of an efficient design of artificial light harvesters \cite{CLCCN,TERNW}.         

The main result of the paper is the following.

\begin{theorem}
Let $w_0 \in L^1(\Rset^3_x \times \Rset^3_{\xi}) \cap L^1(\Rset^3_{\xi}; L^2(\Rset^3_x))$ be a physically admissible initial datum (that is, such that the density matrix operator corresponding to $w_0$ is nonnegative) with $\int_{\Rset^6} |\xi|^2 w_0(x, \xi) \, d(\xi, x) < \infty$. Then, the initial value problem associated with the Unruh-Zurek and the Hu-Paz-Zhang master equations (\ref{UZ}) and (\ref{HPZeq}), respectively, admits a unique global mild solution 
$$
w \in C([0, \infty); L^p(\Rset^3_x \times \Rset^3_{\xi})) \cap C([0, \infty); L^1(\Rset^3_{\xi}; L^2(\Rset^3_x))  \, , \quad 1 \leq p \leq \infty \, .
$$
\label{Th}
\end{theorem}      

The paper is structured as follows: In Section 2 we describe in detail the two main models whose well-posedness is intended to be analyzed, and compute the propagators of the corresponding linear equations as well as establish some of their most relevant properties to our purpose. Section 3 tackles the existence problem of local-in-time mild solutions as well as their main regularity properties. Finally, Section 4 is devoted to investigate the global solvability.  

%%%%%%%%%%%%%%%%%%%%%%%%%%%%%%%%%%%%

%%%%%%%%%%%%%%%%%%%%%%%%%%%%%%%%%%%%%%%%%%%%%%%%%%%%%%%%%%%%%%%%%%%%%%%%%

\section{The Unruh-Zurek and Hu-Paz-Zhang master equations}

%%%%%%%%%%%%%%%%%%%%%%%%%%%%%%%%%%%%%%%%%%%%%%%%%%%%%%%%%%%%%%%%%%%%%%%%%%%%%

The Unruh--Zurek (UZ) master equation was derived in \cite{UZ} to model the interaction of a quantum harmonic oscillator (system) with a scalar field (environment), see also \cite{Zur}. When coupled to the Hartree potential, it  reads as follows:
\begin{equation}
T[w] - \Omega_0^2 (x \cdot \nabla_\xi) w - 2 \gamma \div_{\xi}(\xi w) + \Theta[V]w = c(t) \Delta_{\xi} w - d(t) \div_x(\nabla_\xi w) \, ,
\label{UZ}
\end{equation}
where
\begin{eqnarray*}
c(t) &=& \frac{4\gamma}{\pi} \int_0^\Gamma \theta \coth\left( \frac{\beta \theta}{2} \right) \left( \int_0^t  e^{-\gamma \tau} \cos(\theta \tau) \big( \cos(\Omega \tau) - \frac{\gamma}{\pi} \sin(\Omega \tau) \big)\, d\tau \right) d\theta \, , \\
d(t) &=& \frac{4\gamma}{\Omega \pi} \int_0^\Gamma \theta \coth\left( \frac{\beta \theta}{2} \right) \left( \int_0^t  e^{-\gamma \tau} \cos(\theta \tau) \sin(\Omega \tau) \, d\tau \right) d\theta \, ,
\end{eqnarray*}
and where $\gamma = \frac{\varepsilon^2}{4}$ is the damping coefficient, $\varepsilon$ is the strength of the coupling between the field and the system, and $\Omega = \sqrt{\Omega_0^2 - \frac{\varepsilon^4}{16}}$ denotes the angular frequency of the damped harmonic oscillator, $\Omega_0$ being that of the undamped one. Notice that a units system has been used in which the Planck and Botzmann constants as well as the mass of the oscillator have been normalized to unity.

The Hu--Paz--Zhang (HPZ) master equation was derived for the first time in \cite{HPZ} (see also \cite{HY} for an alternative derivation that corrects a slight deviation in the original coefficients) to model the evolution of the (quasi)distribution function $w(t,x,\xi)$ associated with a Brownian particle (with natural frequency $\Omega$) that interacts linearly with a general environment. When coupled to the Hartree potential, the HPZ equation  reads as follows:
\begin{eqnarray}
T[w] - (\Omega^2 - 2 a(t)) (x \cdot \nabla_\xi) w - 2 b(t) \div_{\xi}(\xi w) + \Theta[V]w = c(t) \Delta_{\xi} w - d(t) \div_x(\nabla_\xi w) \, ,
\label{HPZeq}
\end{eqnarray}
where
\begin{eqnarray}
\label{at}
a(t) &=& \int_0^t  \left( \int_0^\infty I(\theta) \sin(\theta \tau) \, d\theta \right) \cos(\Omega \tau)  \, d\tau \, , \\
\label{bt}
b(t) &=& \frac{1}{\Omega} \int_0^t  \left( \int_0^\infty I(\theta) \sin(\theta \tau) \, d\theta \right) \sin(\Omega \tau)  \, d\tau \, , \\
\label{ct}
c(t) &=& \int_0^t  \left( \int_0^\infty I(\theta) \coth\left( \frac{\beta \theta}{2} \right) \cos(\theta \tau) \, d\theta \right) \cos(\Omega \tau)  \, d\tau \, , \\
\label{dt}
d(t) &=& -\frac{1}{\Omega} \int_0^t  \left( \int_0^\infty I(\theta) \coth\left( \frac{\beta \theta}{2} \right) \cos(\theta \tau) \, d\theta \right) \sin(\Omega \tau)  \, d\tau \, ,
%c(t) &=& \frac{4\gamma}{\pi} \int_0^\Gamma \omega \coth\left( \frac{\beta \omega}{2} \right) \left( \int_0^t  e^{-\gamma \tau} \cos(\omega \tau) \big( \cos(\Omega \tau) - \frac{\gamma}{\pi} \sin(\Omega \tau) \big)\, d\tau \right) d\omega \, , \\
%d(t) &=& \frac{4\gamma}{\Omega \pi} \int_0^\Gamma \omega \coth\left( \frac{\beta \omega}{2} \right) \left( \int_0^t  e^{-\gamma \tau} \cos(\omega \tau) \sin(\Omega \tau) \, d\tau \right) d\omega \, ,
\end{eqnarray}
$\beta = \frac{1}{T}$ is the inverse temperature of the bath, and where $I(\theta)$ denotes the spectral density of the environment. Here, we are intended to deal with an Ohmic environment 
%(that is, $I(\omega) \approx \omega$ in the physical range of frequencies $\omega < \Lambda$, with $\Lambda$ denoting a certain cutoff frequency for which $I(\omega) \to 0$ when $\omega > \Lambda$) 
described by the following spectral density (Drude-Lorentz cutoff):
$$
I(\theta) = \frac{\delta \Gamma^2  \theta}{\Gamma^2 + \theta^2}  \, , \quad \delta > 0 \, .
$$
%that vanishes for values of $\omega$ above the cutoff frequency $\Gamma$.

By a mild solution of the initial value problem associated with Eq. (\ref{UZ}) (resp. (\ref{HPZeq})), subjected to the initial condition $w(0,x, \xi) = w_0(x,\xi)$,  we understand a continuous function $w: [0, T] \rightarrow X$ (the function space $X$ will be specified later on) satisfying the following integral equation:
\begin{eqnarray}
\nonumber
w(t,x,\xi) &=& \int_{\Rset^6} G(t, x,\xi,z,v) w_0(z,v) \, d(z,v) \\ &-& \int_0^t \int_{\Rset^6} G(t-s,x,\xi,z,v)  (\Theta[V]w)(s,z,v) \, d(z,v) \, ds \, ,
\label{mild}
\end{eqnarray}
%with
%$$
%L = \Theta[V]w - 2A(s) (z \cdot \nabla_v) w + \frac{1}{\Omega} B(s) \div_{v}(v w) - 2 C(s) \Delta_{v} w + \frac{2}{\Omega} D(s) \div_z(\nabla_v W) \, ,
%$$
$G(t,x,\xi,z,v)$ denoting the fundamental solution (or propagator) of the linear UZ initial value problem $L_{UZ}[w] = 0, w(0, x, \xi) = w_0(x,\xi)$ (resp. the linear HPZ initial value problem $L_{HPZ}[w] = 0, w(0, x, \xi) = w_0(x,\xi)$), where 
\begin{equation}
\label{linUZ} 
L_{UZ}[w] = T[w] - \Omega_0^2 (x \cdot \nabla_\xi) w - 2 \gamma \div_{\xi}(\xi w)- c(t) \Delta_{\xi} w + d(t) \div_x(\nabla_\xi w) \, , 
\end{equation}
\begin{equation}
\label{linHPZ}
L_{HPZ}[w] = T[w] - (\Omega^2 - 2 a(t)) (x \cdot \nabla_\xi) w - 2 b(t) \div_{\xi}(\xi w) - c(t) \Delta_{\xi} w + d(t) \div_x(\nabla_\xi w) \, .
\end{equation}
Thus,  we may consider the distribution function $w(t,x,\xi)$ to be split into two parts: a linear part, only depending on the initial data $w_0$, and a nonlinear part depending upon the potential $V$ through the pseudo--differential operator $\Theta[V]w$. 

%%%%%%%%%%%%%%%%%%%%%%%%%%

\subsection{The linear UZ problem}

%%%%%%%%%%%%%%%%%%%%%%%%%%%%%%%%%%%%%%%

In order to have a suitable description of (\ref{mild}), our first aim is to compute the fundamental solution $G_{UZ}(t, x,\xi,z,v)$ of the linear problem. Fourier transforming $L_{UZ}[w] = 0$ (cf. \ref{linUZ})) one arrives at
 \begin{equation}
 \widehat{G_0^{UZ}}(t, Y(t), \Pi(t)) = e^{\int_0^t d(s) Y(s) \cdot \Pi(s) \, ds - \int_0^t c(s) |\Pi(s)|^2 \, ds}  \, ,
 \label{FUZ}
 \end{equation}
 where $(Y(t), \Pi(t))$ satisfies the following system of characteristic equations:
 $$
 \left\{
 \begin{array}{ll}
 Y'(t) = \Omega_0^2 \, \Pi(t) \\
\Pi'(t) = 2 \gamma \Pi(t) - Y(t) 
\end{array}
 \, ,
 \right.
 $$
that might be considered to be subjected to the initial conditions $Y(0) = y$ and $\Pi(0) = \eta$, so that $\Pi'(0) = 2 \gamma \eta - y$. Then, by taking a new time derivative in the second equation the characteristic system can be rewritten as
$$
 \left\{
 \begin{array}{ll}
Y'(t) = \Omega_0^2 \, \Pi(t)  \\
\Pi''(t) - 2 \gamma \Pi'(t) + \Omega_0^2 \, \Pi(t) = 0 
\end{array}
 \, ,
 \right.
 $$
 that can be solved explicitely to find
 \begin{eqnarray}
 \label{PiUZ}
 \Pi(t) &=& e^{\frac{\gamma}{2} t} \left( \eta \cos(\Omega t) + \frac{1}{\Omega} \left(  \frac{3}{2} \gamma \eta - y\right)  \sin(\Omega t) \right) \, , \\
 \label{YUZ}
 Y(t) &=& y + \Omega_0^2 \int_0^t \Pi(s) \, ds = y + \Omega_0^2 \left( \eta \alpha(t) + \frac{1}{\Omega} \left(  \frac{3}{2} \gamma \eta - y\right) \beta(t)\right) \, ,
 \end{eqnarray}
 with
 \begin{eqnarray*}
\alpha(t) &=& \frac{2\gamma}{\gamma^2 + 4 \Omega^2} \left(  \frac{2\Omega}{\gamma} e^{\frac{\gamma}{2} t} \sin(\Omega t) + e^{\frac{\gamma}{2} t} \cos(\Omega t) -1 \right) \, , \\
\beta(t) &=&  \frac{2\gamma}{\gamma^2 + 4 \Omega^2} \left(  e^{\frac{\gamma}{2} t} \sin(\Omega t) - \frac{2\Omega}{\gamma}  e^{\frac{\gamma}{2} t} \cos(\Omega t) + \frac{2\Omega}{\gamma} \right) \, . 
 \end{eqnarray*}
 
%thus 
%$$
%Y(t) \equiv y \, , \quad \Pi(t) = \eta e^{2 \gamma t} + \frac{1}{2\gamma} \big( 1 - e^{2\gamma t} \big) y \, .
%$$

%%%%%%%%%%%%%%%%%%%%%%%%%%%%%%%%%%%%%%%%%%%%%%%%
%
%\subsection{The free particle case ($\Omega_0 = 0$)}
%
%%%%%%%%%%%%%%%%%%%%%%%%%%%%%
%
%In this case we have $\Omega = i \gamma$, so that the coefficients $c(t)$ and $d(t)$ can be rewritten as

For the sake of simplicity we restrict ourselves to the free particle model, namely $\Omega_0 = 0$, in which case $\Omega = i \gamma$. Then,
the characteristic curves coincide with those computed for the Markovian Wigner-Fokker-Planck equation (MWFPE) studied in \cite{CLN}:
\begin{equation}
Y(t) \equiv y \, , \quad \Pi(t) = e^{2 \gamma t} \eta - \frac{1}{2\gamma} (e^{2\gamma t} -1) y \, ,
\label{cc}
\end{equation}
while the coefficients $c(t)$ and $d(t)$, that keep constant for the MWFPE, are now reduced to the following time-dependent expressions:
\begin{eqnarray*}
c(t) &=& \frac{4\gamma}{\pi} \int_0^\Gamma \theta \coth\left( \frac{\beta \theta}{2} \right) \left( \int_0^t  e^{-\gamma \tau} \cos(\theta \tau) \big( \cosh(\gamma \tau) - \sinh(\gamma \tau) \big)\, d\tau \right) d\theta  \\ &=& \frac{4\gamma}{\pi} \int_0^\Gamma \theta \coth\left( \frac{\beta \theta}{2} \right) \left( \int_0^t  e^{-2\gamma \tau} \cos(\theta \tau) \, d\tau \right)  d\theta \\ &=& \frac{4\gamma}{\pi} \left( 2 \gamma I_1'(0) - \big(2 \gamma I_1'(t) + I_1''(t) \big) e^{-2 \gamma t} \right) \, ,  \\
d(t) &=& \frac{4}{\pi} \int_0^\Gamma \theta \coth\left( \frac{\beta \theta}{2} \right) \left( \int_0^t  e^{-\gamma \tau} \cos(\theta \tau) \sinh(\gamma \tau) \, d\tau \right) d\theta \\ &=& \frac{2}{\pi} \int_0^\Gamma \theta \coth\left( \frac{\beta \theta}{2} \right) \left( \int_0^t  \cos(\theta \tau) \big( 1- e^{-2\gamma \tau} \big)  \, d\tau \right) d\theta \\ &=& - \frac{1}{2\gamma} c(t) + \frac{2}{\pi} \big(  4\gamma^2 I_1(t) - I_1''(t) \big) \, ,
\end{eqnarray*}
where we denoted
$$
I_1(t) = \int_0^{\Gamma} \frac{1}{4 \gamma^2 + \theta^2} \coth\left( \frac{\beta \theta}{2} \right) \sin(\theta t) \, d\theta \, .
$$
In the following we consider only the case $T=0$ (or equivalently $\beta = \infty$, and thus $\coth\left( \frac{\beta \theta}{2} \right) = 1$) to avoid the logarithmic divergence of $I_1(t)$, as discussed in \cite{UZ}. Also, we assume that the cutoff frequency $\Gamma > 0$ is sufficiently small so that the terms of order $O(\theta^3)$ become negligible. Under these assumptions we can easily observe that $I_1(t)$ grows linearly with time:  
$$
I_1(t) = I_1'(0) t \, ,
$$
for which the truncated Taylor expansion $\sin(\theta t) = \theta t + O(\theta^3)$ has been taken into account. Also, the expressions for $c(t)$ and $d(t)$ become
\begin{eqnarray*}
c(t) &=& \frac{8\gamma^2}{\pi} I_1'(0) \big( 1 - e^{-2\gamma t}\big) \, , \\
d(t) &=& - \frac{1}{2\gamma} c(t) + \frac{8 \gamma^2}{\pi} I_1'(0) t = \frac{4\gamma}{\pi} I_1'(0) \big( e^{-2\gamma t} + 2\gamma t -1 \big)  \, .
\end{eqnarray*}

Now, inserting (\ref{cc}) into Eq. (\ref{FUZ}) yields
 \begin{equation}
 \widehat{G_0^{UZ}}(t, Y(t), \Pi(t)) = e^{- A_{UZ}(t) |y|^2 - B_{UZ}(t) y \cdot \eta - C_{UZ}(t) |\eta|^2}  \, ,
 \label{FUZbis}
 \end{equation}
 with
 \begin{eqnarray*}
  A_{UZ}(t) 
  &=& \frac{1}{4\gamma^2} \int_0^t c(s) (e^{2\gamma s}-1)^2 \, ds + \frac{1}{2\gamma} \int_0^t d(s) (e^{2 \gamma s} -1) \, ds \, , \\  
%  &=& \frac{4\gamma}{\pi} \big( I_1'(0) - 8 \gamma^2 I_2'(0) \big) + \frac{8\gamma^2}{\pi} I_1'(0) t - \frac{4\gamma}{\pi} \big( 4 \gamma I_1(t) + I_1'(t) - 16 \gamma^3 I_2(t) - 8 \gamma^2 I_2'(t)\big) e^{-2 \gamma t} \\ &=& - \frac{4\gamma}{\pi} \big[ I_1'(0) (1-2\gamma t) + \big(  I_1'(t) - 4 \gamma I_2''(t) - 8 \gamma^2 I_2'(t)\big) e^{-2 \gamma t} \big] \, , \\
  B_{UZ}(t) &=& - \frac{1}{\gamma}  \int_0^t c(s) e^{2\gamma s} (e^{2\gamma s}-1) \, ds - \int_0^t d(s) e^{2\gamma s} \, ds  \, , \\ 
%  &=& - \frac{8\gamma}{\pi} \big[ \big( 2 I_1'(0) + 6\gamma I_1(t) + 3 I_1'(t)\big) - \big( 3 I_1'(0) - 2 \gamma I_1(t) + 3 I_1'(t)  \big) e^{2\gamma t} + I_1'(0) e^{4\gamma t} \big] \, , \\
  C_{UZ}(t) &=& \int_0^t c(s) e^{4\gamma s} \, ds  \, .
  %= \frac{8 \gamma^2}{\pi} \big( I_1'(0) - 2 I_1'(t) e^{2\gamma t} + I_1'(0) e^{4\gamma t}\big) \, ,
 \end{eqnarray*}
 As consequence, the change of variables $z=y$, $v = e^{2\gamma t} \eta - \frac{1}{2\gamma} (e^{2\gamma t} -1) z$ leads to
 \begin{equation}
  \widehat{G_0^{UZ}}(t, z, v) = e^{- \widetilde{A_{UZ}}(t) |z|^2 - \widetilde{B_{UZ}}(t) z \cdot v - \widetilde{C_{UZ}}(t) |v|^2}
 \label{hatG}
 \end{equation}
 with
 \begin{eqnarray}
 \label{tildeA}
 \widetilde{A_{UZ}}(t) &=& A_{UZ}(t) + \frac{1}{2\gamma} (1- e^{-2\gamma t}) B_{UZ}(t) + \frac{1}{4\gamma^2} (1- e^{-2\gamma t})^2 C_{UZ}(t) \, , \\
  \label{tildeB}
  \widetilde{B_{UZ}}(t) &=& e^{-2\gamma t} \Big( B_{UZ}(t) + \frac{1}{\gamma} (1-e^{-2\gamma t}) C_{UZ}(t) \Big) \, , \\
  \label{tildeC}
  \widetilde{C_{UZ}}(t) &=& e^{-4 \gamma t} C_{UZ}(t) \, .
 \end{eqnarray}

The fundamental solution is then computed as established in the following 
 \begin{lemma}[Fundamental solution of the linear Unruh-Zurek master equation]
Assume that $\Gamma > 0$ is sufficiently small in the sense stated above. Then, the coefficients $A_{UZ}(t)$, $B_{UZ}(t)$, and $C_{UZ}(t)$ of $\widehat{G_0^{UZ}}$ are given by
 \begin{eqnarray}
 \label{A-UZ}
 A_{UZ}(t) &=& \frac{1}{2 \gamma \pi} I_1'(0) \Big(  4 e^{-2\gamma t} + \left(4\gamma t - 10\right) e^{2\gamma t} + e^{4\gamma t} + 20\gamma t - 4 \gamma^2 t^2+ 5 \Big) \, , \\
 \label{B-UZ}
 B_{UZ}(t) &=&  -\frac{2}{\pi} I_1'(0) \Big(  \left(2\gamma t -6\right) e^{2\gamma t} + e^{4\gamma t} + 6\gamma t + 5 \Big) \, ,\\
 \label{C-UZ}
 C_{UZ}(t) &=& \frac{2\gamma}{\pi} I_1'(0) (e^{2\gamma t}-1)^2 \, .
 \end{eqnarray}
 As consequence, 
\begin{eqnarray}
\label{tildes1}
\widetilde{A_{UZ}}(t) &=&  \frac{2}{\pi} I_1'(0) F_1(t) \, ,  \\ \label{tildes2} \widetilde{B_{UZ}}(t) &=& -\frac{2}{\pi} I_1'(0) F_2(t) \, ,  \\ \label{tildes3 }\widetilde{C_{UZ}}(t) &=& \frac{2\gamma}{\pi} I_1'(0) e^{-4\gamma t} (e^{2\gamma t}-1)^2  \, ,
\end{eqnarray}
 with
 \begin{eqnarray*}
F_1(t) &=&  \frac{1}{4\gamma} \Big( 10 e^{-2\gamma t} + e^{-4\gamma t} - 11 + 4\gamma t (3e^{-2\gamma t} + 3 - \gamma t) \Big)\, , \\
F_2(t) &=& 2 e^{-2\gamma t} + e^{-4\gamma t} - 3 + 2\gamma t (3e^{-2\gamma t} + 1) \, .
 \end{eqnarray*}
 Then, the fundamental solution of the linear Unruh-Zurek master equation is given by 
 \begin{equation}
 G_{UZ}(t,x,\xi,z,v) = G_0^{UZ}\left(t, x-z-\frac{1}{2\gamma} (1-e^{-2\gamma t})v, \xi - e^{-2\gamma t} v\right) \, ,
 \label{G}
 \end{equation}
 where
 \begin{equation}
 G_0^{UZ}(t,x,\xi) = d_{UZ}(t) \, e^{-a_{UZ}(t)|x|^2 + b_{UZ}(t) x \cdot \xi - c_{UZ}(t) |\xi|^2} \, , 
 \label{G0}
 \end{equation}
 with
 \begin{equation}
 d_{UZ}(t) = \frac{1}{(2\pi)^3} \big( 4 \widetilde{A_{UZ}}(t) \widetilde{C_{UZ}}(t) - \widetilde{B_{UZ}}^2 \big)^{-\frac{3}{2}}
 \label{d}
 \end{equation}
 and
 \begin{eqnarray}
 \label{a}
 a_{UZ}(t) &=& 4\pi^2 d_{UZ}(t)^{\frac{2}{3}} \widetilde{C_{UZ}}(t)  \, ,  \\  \label{b} b_{UZ}(t) &=& 4\pi^2 d_{UZ}(t)^{\frac{2}{3}} \widetilde{B_{UZ}}(t) \, , \\ \label{c} c_{UZ}(t) &=& 4\pi^2 d_{UZ}(t)^{\frac{2}{3}}  \widetilde{A_{UZ}}(t) \, .
 \label{abc}
 \end{eqnarray}
 \end{lemma}
 
\noindent {\bf{Proof.}} 
We have
 \begin{eqnarray*}
  A_{UZ}(t) &=& \frac{2}{\pi} I_1'(0) \left\{ \int_0^t  (1 - e^{-2\gamma s}) (1 - e^{2\gamma s})^2 \, ds +  \int_0^t \big( e^{-2\gamma s} + 2\gamma s -1 \big) (e^{2 \gamma s} -1) \, ds \right\} \\
&=&  \frac{2}{\pi} I_1'(0) \left\{ \int_0^t  (5 - 4 e^{2\gamma s} + e^{4\gamma s} - 2e^{-2\gamma s}  + 2\gamma s e^{2\gamma s} -2\gamma s ) \, ds \right\} \\ &=& \frac{2}{\pi} I_1'(0) \left( 5t - \gamma t^2 + \frac{1}{\gamma} e^{-2\gamma t} + \left(t - \frac{5}{2\gamma}\right) e^{2\gamma t} + \frac{1}{4\gamma} e^{4\gamma t} + \frac{5}{4\gamma}  \right) \, ,
%  - \frac{2}{\gamma \pi} \left( \frac{5}{2} I_1''(0) + 10 \gamma^2 I_1(0) + 2\gamma I_1''(0) t - 2\gamma^2 (1+16\gamma^2)I_1(t) + 3 \gamma I_1'(t) + 2 I_1''(t) \right) \\ &-& \frac{1}{2\gamma \pi}  I_1''(0) e^{4\gamma t} + \frac{1}{\gamma \pi}  \big( 3 I_1''(0) + 2 I_1''(t) - 2 \gamma I_1'(t) \big) e^{2 \gamma t}  \\ &+& \frac{2}{\gamma \pi} \big( 8 \gamma^2 I_1(t) + 4 \gamma I_1'(t) + I_1''(t) - 32 \gamma^4 I_2(t) - 16 \gamma^3 I_2'(t)\big) e^{-2 \gamma t} \, .
 \end{eqnarray*}
 that straightforwardly leads to (\ref{A}). We now compute 
  \begin{eqnarray*}
  B_{UZ}(t) &=& -\frac{4 \gamma}{\pi} I_1'(0) \left\{ 2 \int_0^t (e^{2\gamma s} -1)^2 \, ds + \int_0^t \big(1 + 2\gamma s e^{2\gamma s} - e^{2\gamma s}\big) \, ds \right\} 
  \\ &=& -\frac{4\gamma}{\pi} I_1'(0) \left( \frac{1}{2\gamma} e^{4\gamma t} + \left(t -\frac{3}{\gamma}\right) e^{2\gamma t} + 3 t + \frac{5}{2\gamma} \right) \, ,
  %&=&  \frac{2}{\pi} \big( 6\gamma I_1'(t) + 3 I_1''(t) + 2 I_1''(0)\big) +  \frac{2}{\pi}  I_1''(0) e^{4\gamma t}  \\ &-& \frac{2}{\pi} \big(3 I_1''(0) + 3 I_1''(t) - 2 \gamma I_1'(t) \big) e^{2\gamma t} \, .
 \end{eqnarray*}
 that is equivalent to (\ref{B}). Finally, we have
  \begin{eqnarray*}
  C_{UZ}(t) &=&  \frac{2\gamma}{\pi} I_1'(0) (e^{2\gamma t}-1)^2 \, .
 \end{eqnarray*}
 % where 
% $$
% I_2(t) = \int_0^{\Gamma} \frac{1}{(4 \gamma^2 + \omega^2)^2} \coth\left( \frac{\beta \omega}{2} \right) \sin(\omega t) \, d\omega \, , 
% $$
% so that the following relations hold: $I_1(t) = 4 \gamma^2 I_2(t) - I_2''(t)$, $I_1'(0) = 4\gamma^2 I_2'(0)$.
Then, the coefficients with tilde are directly obtained from their definitions (\ref{tildeA})-(\ref{tildeC}). 

Now, by applying the inverse Fourier transform to $ \widehat{G_0^{UZ}}$ (cf. (\ref{hatG})) we find
 $$
 G_0^{UZ}(t,x,\xi) = d_{UZ}(t) \, e^{-a_{UZ}(t)|x|^2 + b_{UZ}(t) x \cdot \xi - c_{UZ}(t) |\xi|^2} \, , 
 $$
 with the coefficients given as in (\ref{d})-(\ref{c}).
 Clearly $\widetilde{C_{UZ}} > 0$, and therefore $ a_{UZ} > 0$. Also, since  $I_1'(0) =  \frac{1}{2} \ln\left( 1 + \frac{\Gamma^2}{4\gamma^2}\right) > 0$, it suffices to show that $F_1>0$ in order to conclude that $\widetilde{A_{UZ}} > 0$. Indeed, $F_1(0) = 0$ and 
 $$
 F_1\left(\frac{1}{2\gamma}\right) = \frac{1}{4\gamma} \left( \frac{16}{e} + \frac{1}{e^2} -5\right) > 0 
 $$
 for arbitrary $\gamma > 0$. 
 Hence, $F_1(t) > 0$ for all $t>0$. Finally, we have to check that $4 \widetilde{A_{UZ}}(t) \widetilde{C_{UZ}}(t) - \widetilde{B_{UZ}}(t)^2 > 0$ in order to give sense to $ d_{UZ}(t)$. A simple computation leads us to the relation $4 \widetilde{A_{UZ}}(t) \widetilde{C_{UZ}}(t) - \widetilde{B_{UZ}}(t)^2 = e^{-4 \gamma t} \big( 4 A_{UZ}(t) C_{UZ}(t) - B_{UZ}(t)^2 \big)$. Now, inserting the expressions (\ref{A-UZ})-(\ref{C-UZ}) for the coefficients without tilde and rearranging terms properly we obtain
$$
 4 \widetilde{A_{UZ}}(t) \widetilde{C_{UZ}}(t) - \widetilde{B_{UZ}}(t)^2 = \frac{16}{\pi^2} I_1'(0)^2 e^{-4 \gamma t} H_{UZ}(t) \, ,
$$
 with
$$
 H_{UZ}(t) = (6 \gamma t - 2 \gamma^2 t^2 -5) e^{4\gamma t} + (4 \gamma t - 4 \gamma^2 t^2 +11) e^{2\gamma t} + e^{-2 \gamma t} - 10\gamma t (1 + \gamma t) - 7 \, .
$$
Again, $H_{UZ}(0) = 0$ and 
$$
 H_{UZ}\left(\frac{1}{2\gamma}\right) = -\frac{1}{2} (5e^{2} + 29) + 12 e + \frac{1}{e} > 0 \, ,
$$
thus $H_{UZ}(t) > 0$ for all $t>0$, and consequently $d_{UZ}(t)$ is shown to be well defined.
 
  \hfill 
 $\square$
% The consistency of the fundamental solution is checked through the following 
% \begin{lemma}
% $d_{UZ}(t) > 0$ for all $t > 0$.
% \label{consist}
% \end{lemma}
% 
%\noindent  {\bf{Proof.}} 
%
% \hfill 
% $\square$

We now establish some integrability and regularity properties of the fundamental solution that will be helpful later on. In what follows, $C(p,q)$ will denote a positive constant depending on the integrability indices $p$ and $q$, and 
$$
\| f\|_{L^{q,p}} = \| f \|_{L^q(\Rset^3_\xi; L^p(\Rset^3_x))} = \left( \int_{\Rset^3} \left( \int_{\Rset^3} | f(x,\xi) |^p \, dx \right)^{\frac{q}{p}} \, d\xi \right)^{\frac{1}{q}} \, .
$$
%$$
%\| f(t) \|_{L^{q,p}} = \Big(  \int_{\Rset^3} \Big(\int_{\Rset^3} |f(t,x,\xi)|^p \, dx \Big)^{\frac{q}{p}} \, d\xi \Big)^{\frac{1}{q}} \, .
%$$

\begin{lemma}[Properties of $G_{UZ}$]

The fundamental solution of the Unruh-Zurek master equation, given by (\ref{G})-(\ref{c}), satisfies the following properties for all times $t > 0$:  
\begin{itemize}
\item[(i)] $\int_{\Rset^6} G_{UZ}(t,x,\xi,z,v) \, d(\xi,x) = 1$ for all $(z,v) \in \Rset^6$.
\item[(ii)] $\|G_0^{UZ} (t)\|_{L^{q,p}} \leq C(q,p) \, a_{UZ}(t)^{\frac{3}{2}(\frac{1}{q} - \frac{1}{p})} \, d_{UZ}(t)^{1-\frac{1}{q}}$ for all $1 \leq q , p < \infty$.
\item[(iii)] Given $1 \leq q, p < \infty$, the equality 
$$
\big\| |\xi|^{\alpha} G_0^{UZ}(t) \big\|_{L^{q,p}} = C(p,q) \, a_{UZ}(t)^{\frac{\alpha}{2} - \frac{3}{2}(\frac{1}{p} - \frac{1}{q})} \, d_{UZ}(t)^{1-\frac{1}{q}-\frac{\alpha}{3}} \, .
$$
is fulfilled.
%\item[(iv)] $\| \nabla_x G_0^{UZ} (t)\|_{L^{q,p}} \leq C(q,p) \, a_{UZ}(t) \, d_{UZ}(t)^{\frac{2}{3} - \frac{1}{p}} \, c_{UZ}(t)^{\frac{1}{2} \left( 1+\frac{3}{p} - \frac{3}{q}\right)}$.
%\begin{eqnarray*}
%\| \nabla_{(x,\xi)}G_0^{UZ} (t)\|_{L^p(\Rset^3_x;L^q(\Rset^3_\xi))} &\leq& C(q,p) \Big(  \big(2a_{UZ}(t) + b_{UZ}(t)\big) d_{UZ}(t)^{\frac{2}{3} - \frac{1}{p}} c_{UZ}(t)^{\frac{1}{2} \left( 1+\frac{3}{p} - \frac{3}{q}\right)}  \\ && + \,  \big(2c_{UZ}(t) + b_{UZ}(t)\big) d_{UZ}(t)^{\frac{2}{3} - \frac{1}{q}} a_{UZ}(t)^{\frac{1}{2} \left( 1+ \frac{3}{q} - \frac{3}{p}\right)}\Big) \, .
%\end{eqnarray*}
\end{itemize}
%Here, we denoted $C(p,q)$ a positive constant that only depends on $p$ and $q$.
\label{lema2}
\end{lemma}

\noindent {\bf{Proof.}} Denote $X(t,x,z,v) = x-z-\frac{1}{2\gamma} (1-e^{-2\gamma t})v$. The first assertion follows from
\begin{eqnarray*}
&&\int_{\Rset^6} G_{UZ}(t,x,\xi,z,v) \, d(\xi,x) = \int_{\Rset^6} G_0^{UZ} \left(t, X(t,x,z,v), \xi - e^{-2\gamma t} v\right) \, d(\xi,x) \\ && = d_{UZ}(t) \, e^{- c_{UZ}(t) e^{-4\gamma t} |v|^2} \int_{\Rset^3} e^{-a_{UZ}(t) \left| X(t,x,z,v) \right|^2}  e^{-b_{UZ}(t) e^{-2\gamma t} X(t,x,z,v) \cdot v} J(t,x,z,v) \, dx \, ,
\end{eqnarray*}
where
\begin{eqnarray*}
J(t,x,z,v) &=& \int_{\Rset^3} e^{- c_{UZ}(t) |\xi|^2} e^{ j(t,x,z,v) \cdot \xi} \, d\xi = \pi^{\frac{3}{2}} c_{UZ}(t)^{- \frac{3}{2}} \, e^{\frac{|j(t,x,z,v)|^2}{4 c_{UZ}(t)}}
\\ &=& \pi^{\frac{3}{2}} c_{UZ}(t)^{- \frac{3}{2}} \, e^{c_{UZ}(t) e^{-4\gamma t} |v|^2} e^{b_{UZ}(t) e^{-2\gamma t} X(t,x,z,v) \cdot v} e^{\frac{b_{UZ}(t)^2}{4 c_{UZ}(t)} |X(t,x,z,v)|^2}
\end{eqnarray*}
with
$$
j(t,x,z,v) =  2 c_{UZ}(t) e^{-2\gamma t} v + b_{UZ}(t) X(t,x,z,v) \, .
$$
Then,
\begin{eqnarray*}
&& \int_{\Rset^3} e^{-a_{UZ}(t) \left| X(t,x,z,v) \right|^2}  e^{-b_{UZ}(t) e^{-2\gamma t} X(t,x,z,v) \cdot v} J(t,x,z,v) \, dx\\ &&= \pi^{\frac{3}{2}} c_{UZ}(t)^{- \frac{3}{2}} \, e^{c_{UZ}(t) e^{-4\gamma t} |v|^2} \int_{\Rset^3} e^{-\left(a_{UZ}(t) - \frac{b_{UZ}(t)^2}{4 c_{UZ}(t)} \right) |X(t,x,z,v)|^2} \, dx \\ &&= \pi^{\frac{3}{2}} c_{UZ}(t)^{- \frac{3}{2}} \, e^{c_{UZ}(t) e^{-4\gamma t} |v|^2} \int_{\Rset^3} e^{-\pi^2 \frac{d_{UZ}(t)^{\frac{2}{3}}}{c_{UZ}(t)} |x|^2} \, dx \\ &&= \frac{1}{d_{UZ}(t)} \, e^{c_{UZ}(t) e^{-4\gamma t} |v|^2} \, ,
\end{eqnarray*}
for which we used the relations (\ref{d})-(\ref{c}). 

Assertions (ii) and (iii) follow from straightforward computations involving Gaussian integrals (similar to those of (i)) and taking into account formulae (\ref{d})-(\ref{c}) again.

%Finally, we observe that 
%$$
%\| \nabla_x G_0^{UZ} (t)\|_{L^{q,p}} \leq 2a_{UZ}(t) \| x G_0^{UZ} (t)\|_{L^{q,p}} \leq 2a_{UZ}(t) \| x G_0^{UZ} (t)\|_{L^p(\Rset^3_x; L^q(\Rset^3_\xi))} 
%$$
%by virtue of Minkowski's inequality for integrals, which allows to reverse the order of the norms. The computation of $\| x G_0^{UZ} (t)\|_{L^{q,p}}$ follows the same steps that in (iii), but now changing the integration order is not required.   
%Then,
%\begin{eqnarray*}
%\big\| x G_0^{UZ} (t) \big\|_{L^p(\Rset^3_x; L^q(\Rset^3_\xi))} &=& d_{UZ}(t) \left(  \int_{\Rset^3_x} |x|^p e^{-p \, a_{UZ}(t) |x|^2} \left(  \int_{\Rset^3_\xi} e^{-q \, c_{UZ}(t) |\xi|^2 + q \, b_{UZ}(t) \xi \cdot x} \, d\xi \right)^{\frac{p}{q}} \, dx \right)^{\frac{1}{p}} \\ &=& C(p,q) d_{UZ}(t) c_{UZ}(t)^{-\frac{3}{2q}} \left(  \int_{\Rset^3_x} |x|^p e^{-p \pi^2 \, \frac{d_{UZ}(t)^{\frac{2}{3}}}{c_{UZ}(t)} |x|^2} \, dx \right)^{\frac{1}{p}} \\ &=& C(p,q) d_{UZ}(t)^{\frac{2}{3} - \frac{1}{p}}  c_{UZ}(t)^{ \frac{1}{2} + \frac{3}{2p} -\frac{3}{2q} } \, .
%\end{eqnarray*}

This concludes the proof.
%Finally, we have
%\begin{eqnarray*}
%&&\int_{\Rset^6} G(t,x,\xi,z,v) \, d(\xi,x) = \big( 4 \widetilde{A_{UZ}}(t) \widetilde{C_{UZ}}(t) - \widetilde{B_{UZ}}^2 \big)^{-\frac{3}{2}}  (4 a_{UZ}(t) c_{UZ}(t) - b_{UZ}(t)^2)^{-\frac{3}{2}} \\ &&\times e^{ -\frac{1}{4\gamma^2} ( a_{UZ}(t) - \delta_1(t)) (1-e^{-2\gamma t})^2 |v|^2} e^{\frac{1}{4\gamma^2} (1-e^{-2\gamma t})^2(a_{UZ}(t)  -  \delta_1(t)) z \cdot v} 
%\end{eqnarray*}

\hfill 
 $\square$
 
 The following result will be useful for future estimates, too.
 
 \begin{lemma}[Rates of time growth/decay near $t=0$ of the coefficients of $G_{UZ}$] Let $a_{UZ}(t)$, $b_{UZ}(t)$, $c_{UZ}(t)$ and $d_{UZ}(t)$ be given by (\ref{d})-(\ref{c}) and denote 
 $$
 D_{UZ}(t) = 4 A_{UZ}(t) C_{UZ}(t) - B_{UZ}(t)^2 \, . 
 $$
 Then, their behavior near $t=0$ is as follows:   
 $$
D_{UZ}(t) \sim t^6 \, , \quad  a_{UZ}(t)  \sim t^{-4} \, , \quad b_{UZ}(t) \sim t^{-3} \, , \quad c_{UZ}(t) \sim t^{-2} \, , \quad d_{UZ}(t) \sim t^{-9} \, .
 $$
 \label{lema3}
 \end{lemma}
 
 \noindent {\bf{Proof.}} We have
 \begin{eqnarray*}
&& A_{UZ}(0) = A_{UZ}'(0) = A_{UZ}''(0) = A_{UZ}'''(0) = 0, \quad A_{UZ}^{(iv)}(0) = \frac{144 \gamma^4}{\pi} I_1'(0) \, , \\
&& B_{UZ}(0) = B_{UZ}'(0) = B_{UZ}''(0) = 0, \quad B_{UZ}'''(0) = - \frac{80 \gamma^4}{\pi} I_1'(0) \, , \\
&& C_{UZ}(0) = C_{UZ}'(0) = 0, \quad C_{UZ}''(0) = \frac{16 \gamma^4}{\pi} I_1'(0) \, , 
\end{eqnarray*}
thus
$$
D_{UZ}(0) =  D_{UZ}'(0) = \cdots = D_{UZ}^{(v)}(0) = 0, \quad D_{UZ}^{(vi)}(0) = 10 \left(\frac{32 \gamma^4}{\pi} I_1'(0)\right)^2 \, ,
$$
from which we may conclude that $D_{UZ}(t)$ goes as $t^6$ in the vicinity of $t=0$.

Also,
$$
\left( \frac{1}{a_{UZ}}\right)(t) = \frac{4 \widetilde{A_{UZ}}(t) \widetilde{C_{UZ}}(t) - \widetilde{B_{UZ}}(t)^2}{\widetilde{C_{UZ}}(t)} = e^{8 \gamma t} \, \frac{D_{UZ}(t)}{C_{UZ}(t)} \, ,
$$
that vanishes at $t=0$. As a matter of fact, the first nonvanishing derivative of $\frac{1}{a_{UZ}}$ at $t=0$ is that of fourth order, $\left( \frac{1}{a_{UZ}}\right)^{(iv)}(0) = \frac{128 \gamma^4}{3 \pi} I_1'(0)$, so that $\frac{1}{a_{UZ}}$ goes as $t^4$ at $t=0$. 
Analogously,
$$
\left( \frac{1}{b_{UZ}}\right)(t) =  \frac{4 \widetilde{A_{UZ}}(t) \widetilde{C_{UZ}}(t) - \widetilde{B_{UZ}}(t)^2}{\widetilde{B_{UZ}}(t)} = \frac{e^{6 \gamma t}  D_{UZ}(t)}{B_{UZ}(t) + \frac{1}{\gamma} (1-e^{-2\gamma t}) C_{UZ}(t)}  \, ,
$$
that leads to $\left( \frac{1}{b_{UZ}}\right)(0) = \left( \frac{1}{b_{UZ}}\right)'(0) =  \left( \frac{1}{b_{UZ}}\right)''(0) = 0$ and $ \left( \frac{1}{b_{UZ}}\right)'''(0) =  \frac{32 \gamma^4}{\pi} I_1'(0)$, as well as
\begin{eqnarray*}
\left( \frac{1}{c_{UZ}}\right)(t) &=& \frac{4 \widetilde{A_{UZ}}(t) \widetilde{C_{UZ}}(t) - \widetilde{B_{UZ}}(t)^2}{\widetilde{A_{UZ}}(t)} \\ &=& \frac{e^{4 \gamma t}  D_{UZ}(t)}{A_{UZ}(t) + \frac{1}{2\gamma} (1- e^{-2\gamma t}) B_{UZ}(t) + \frac{1}{4\gamma^2} (1- e^{-2\gamma t})^2 C_{UZ}(t)} \, ,
\end{eqnarray*}
whose evaluations at $t=0$ give $\left( \frac{1}{c_{UZ}}\right)(0) =  \left( \frac{1}{c_{UZ}}\right)'(0) =  0$, $\left( \frac{1}{c_{UZ}}\right)''(0) = \frac{128 \gamma^4}{3 \pi} I_1'(0)$.
Finally,
$$
\left( \frac{1}{d_{UZ}}\right)(t) = (2\pi)^3 e^{-6\gamma t} D_{UZ}(t)^{\frac{3}{2}} \, ,
$$
so that $\left( \frac{1}{d_{UZ}}\right)(0) = \left( \frac{1}{d_{UZ}}\right)'(0) = \cdots =  \left( \frac{1}{d_{UZ}}\right)^{(viii)}(0) = 0$ and 
$$
\left( \frac{1}{d_{UZ}}\right)^{(ix)}(0) = 420 \sqrt{2} \big(64 \gamma^4 I_1'(0)\big)^3 \, .
$$
Now we are done with the proof.

 \hfill 
 $\square$
 
 %%%%%%%%%%%%%%%%%%%%%%%%%%

\subsection{The linear HPZ problem}

%%%%%%%%%%%%%%%%%%%%%%%%%%%%%%%%%%%%%%%

Fourier transforming $L_{HPZ}[w]=0$ (cf. (\ref{linHPZ})) one arrives at Eq. (\ref{FUZ}),
%\begin{equation}
%\widehat{G_0}(t, Y(t), \Pi(t)) = e^{\int_0^t d(s) Y(s) \cdot \Pi(s) \, ds - \int_0^t c(s) |\Pi(s)|^2 \, ds}  \, ,
%\label{FHPZ}
%\end{equation}
where now $(Y(t), \Pi(t))$ satisfies the following system of characteristic equations:
\begin{equation}
\left\{
\begin{array}{ll}
Y'(t) = (\Omega^2 - 2a(t))  \, \Pi(t) \\
\Pi'(t) = 2 b(t) \Pi(t) - Y(t) 
\end{array}
\, ,
\right.
\label{cs}
\end{equation}
subjected to the initial conditions $Y(0) = y$ and $\Pi(0) = \eta$, so that $\Pi'(0) = - Y(0)= - y$. In the sequel we restrict ourselves to the case of sufficiently small frequencies $\Omega$ and $\Gamma$ (terms of order $O(\Omega^4)$ and $O(\Gamma^4)$ are assumed to be negligible), so that crossed products with the form $\Gamma^{i} \Omega^{j}$ are assumed negligible for $i+j \geq 4$, in order that the propagator can be explicitly calculated. In this approach we find
\begin{eqnarray*}
a(t) &=& \delta \Gamma^2   \int_0^t  \left( \int_0^\infty \frac{\theta}{\Gamma^2 + \theta^2}  \sin(\theta \tau) \, d\theta \right) \cos(\Omega \tau) \, d\tau = \frac{ \delta \pi \Gamma^2}{2}  \int_0^t  e^{-\Gamma \tau} \cos(\Omega \tau) \, d\tau \\ &=& \frac{\delta \pi}{2} \left( \frac{\Gamma^2}{\Gamma^2 + \Omega^2} \right) \Big(  \Gamma + e^{-\Gamma t} \big( \Omega \sin(\Omega t) - \Gamma \cos(\Omega t) \big)\Big) \\ &=& \frac{\delta \pi \Gamma^2}{4}  t (2 - \Gamma t) + {\hbox{\rm{higher order terms.}}}
\end{eqnarray*}
%where we used the truncated expansions $\sin(\Gamma t) = \Gamma t + \frac{\Gamma^3}{6} t^3 + O(\Gamma^5)$, $\cos(\Gamma t) = 1 - \frac{\Gamma^2}{2} t^2 + O(\Gamma^4)$, $\sin(\Omega t) = \Omega t + O(\Omega^3)$, and $\cos(\Omega t) = 1 - \frac{\Omega^2}{2} t^2 + O(\Omega^4)$. 
Similarly,
\begin{eqnarray*}
b(t) &=& \frac{\delta \Gamma^2}{\Omega} \int_0^t  \left( \int_0^\infty \frac{\theta}{\Gamma^2 + \theta^2} \sin(\theta \tau) \, d\theta \right) \sin(\Omega \tau) \, d\tau = \frac{ \delta \pi \Gamma^2}{2 \Omega}  \int_0^t  e^{-\Gamma \tau} \sin(\Omega \tau) \, d\tau  \\ &=& \frac{\delta \pi \Gamma^2}{12}  t^2 (3 - 2 \Gamma t) + {\hbox{\rm{higher order terms.}}} 
\end{eqnarray*}
Then, by taking a new time derivative in the second equation of (\ref{cs}) the characteristic system can be rewritten as
$$
\left\{
\begin{array}{ll}
Y'(t) =\left(\Omega^2 - \frac{\delta \pi \Gamma^2}{2}  t (2 - \Gamma t) \right)  \, \Pi(t)  \\
\Pi''(t) - \frac{\delta \pi \Gamma^2}{6} t^2 (3 - 2 \Gamma t) \Pi'(t) + \left(\Omega^2 - \frac{\delta \pi \Gamma^2}{2} t (4 - 3 \Gamma t) \right) \Pi(t) = 0 
\end{array}
\, .
\right.
$$
%where we used that $b'(t) = \frac{\delta \Lambda^4}{\Omega(4 + \Lambda^4)} \sin(\Omega t) = a(t)$, as follows from a straightforward computation. 

%In the sequel we shall restrict ourselves to the free particle case $\Omega = 0$, in order that the propagator can be explicitly calculated. Then, we have 
%\begin{eqnarray*}
%a(t) &=& \delta \int_0^t  \left( \int_0^\infty \omega e^{-\frac{\omega^2}{\Lambda^2}} \sin(\omega \tau) \, d\omega \right) \, d\tau = -\frac{\delta \Lambda^2}{2} \int_0^t  \left( \int_0^\infty \frac{d}{d\omega} \big(e^{-\frac{\omega^2}{\Lambda^2}}\big) \sin(\omega \tau) \, d\omega \right) \, d\tau \\ &=& - \frac{\delta \Lambda^4}{4} \int_0^t  \left( \int_0^\infty \frac{d}{d\omega} \big(e^{-\frac{\omega^2}{\Lambda^2}}\big) \cos(\omega \tau) \, d\omega \right) \, d\tau = \frac{\Lambda^4}{4} (\delta t - a(t)) \, , 
%\end{eqnarray*}
%so that 
%$$
%a(t) = \lambda_0 t \, , \quad {\hbox{\rm{ with }}} \ \lambda_0 = \frac{\delta \Lambda^4}{4 + \Lambda^4} \, .
%$$
%Thus, the characteristic system is reduced to
%$$
%\left\{
%\begin{array}{ll}
%Y'(t) = - 2 \lambda_0 t  \, \Pi(t)  \\
%\Pi''(t) - f(t) \Pi'(t) - 4 \lambda_0 t \Pi(t)= 0 
%\end{array}
%\, ,
%\right.
%$$
%with
%$$
%f(t) = 2 \lim_{\Omega \to 0} b(t) = \lambda_0 t^2 \, .
%$$
%Up to first order in $\lambda_0$ (that is, when $\delta$ is assumed sufficiently small so that its powers of order $\geq 2$ can be neglected), 
The (unique) solution to the initial value problem associated with the second order equation for $\Pi(t)$ is given by
\begin{equation}
\Pi(t) = f(t) \eta + g(t) y \, , 
\label{pi}
\end{equation}
with
$$
f(t) = 1 - \frac{\Omega^2}{2} t^2 + \frac{\delta \pi \Gamma^2}{3}  t^3 - \frac{\delta \pi \Gamma^3}{8}  t^4 \, , \quad g(t) = -t + \frac{\Omega^2}{6} t^3 - \frac{5 \delta \pi \Gamma^2}{24} t^4 
+  \frac{11 \delta \pi \Gamma^3}{120} t^5 \, . 
$$
Hence,
\begin{equation}
Y(t) = \widetilde{f}(t) \eta + \widetilde{g}(t) y \, ,
\label{Y}
\end{equation}
where
$$
\widetilde{f}(t) = \Omega^2 t - \frac{\delta \pi \Gamma^2}{2} t^2 +  \frac{\delta \pi \Gamma^3}{6} t^3 \, , \quad \widetilde{g}(t) = f(t) \, .
$$

Under our approximation and for the case in which the temperature of the bath is large enough (namely, terms of order $\beta^3$ are negligible), so that $\coth\left( \frac{\beta \theta}{2} \right) =\frac{2}{\beta \theta} + \frac{\beta \theta}{6}$, the coefficients $c(t)$ and $d(t)$ defined by (\ref{ct})-(\ref{dt}) are reduced to the following time-dependent expressions:
\begin{eqnarray*}
c(t) &=& \frac{2\delta \Gamma^2}{\beta} \int_0^t  \frac{d}{d\tau} \left( \int_0^\infty \frac{1}{\theta(\Gamma^2 + \theta^2)} \sin(\theta \tau) \, d\theta \right) \cos(\Omega \tau) \, d\tau \\ &&  + \, \frac{\beta \delta \Gamma^2}{6} \int_0^t  \frac{d}{d\tau} \left( \int_0^\infty \frac{\theta}{\Gamma^2 + \theta^2} \sin(\theta \tau) \, d\theta \right) \cos(\Omega \tau) \, d\tau \\ &=& \frac{\delta \pi \Gamma}{\beta} \int_0^t  e^{-\Gamma \tau} \cos(\Omega \tau) \, d\tau = \frac{\delta \pi \Gamma}{\beta(\Gamma^2 + \Omega^2)} \Big( \Gamma + e^{-\Gamma t} \big( \Omega \sin(\Omega t) - \Gamma \cos(\Omega \tau) \big) \Big) \\ &=& \frac{\delta \pi \Gamma}{6\beta} t \Big( 6- 3\Gamma t + (\Gamma^2 - \Omega^2) t^2\Big) \, , \\
d(t) &=& - \frac{2\delta \Gamma^2}{\beta \Omega} \int_0^t  \frac{d}{d\tau} \left( \int_0^\infty \frac{1}{\theta(\Gamma^2 + \theta^2)} \sin(\theta \tau) \, d\theta \right) \sin(\Omega \tau) \, d\tau \\ &=& - \frac{\delta \pi \Gamma }{\beta \Omega (\Gamma^2 + \Omega^2)} \Big( \Omega - e^{-\Gamma t} \big( \Gamma \sin(\Omega t) + \Omega \cos(\Omega \tau) \big) \Big) \\ &=& - \frac{\delta \pi \Gamma}{24 \beta} t^2 \Big( 12- 8\Gamma t + (3\Gamma^2 - \Omega^2) t^2\Big) \, .
\end{eqnarray*}
%with $\lambda_0 = \delta \left( \Gamma^2 - \frac{\Omega^2}{2} \right)$, $\mu_0 = \frac{\delta \Gamma^2 \Omega^2}{6}$.
%with $\mu_0 = \delta \left( \Gamma^2 - \frac{\Omega^2}{2} \right)$. 
Now, inserting (\ref{pi}) and (\ref{Y}) into Eq. (\ref{FUZ}) yields
\begin{equation}
\widehat{G_0^{HPZ}}(t, Y(t), \Pi(t)) = e^{- A_{HPZ}(t) |y|^2 - B_{HPZ}(t) y \cdot \eta - C_{HPZ}(t) |\eta|^2}  \, ,
 \label{FHPZbis}
\end{equation}
with
\begin{eqnarray}
\nonumber
&&A_{HPZ}(t) =  \int_0^t c(s) g(s)^2 \, ds -  \int_0^t d(s) g(s) \widetilde{g}(s) \, ds \\  \nonumber && = \frac{\delta \pi \Gamma}{24 \beta} \left\{ 12 \int_0^t s g(s) \big(s f(s) + 2 g(s)\big) \, ds - 4 \Gamma \int_0^t s^2 g(s) \big(2s f(s) + 3 g(s)\big) \, ds \right. \\ \label{A-HPZ} && \left.  + \Gamma^2  \int_0^t s^3 g(s) \big(3s f(s) + 4 g(s)\big) \, ds - \Omega^2 \int_0^t s^3 g(s) \big(s f(s) + 4 g(s)\big) \, ds  \right\} \, , \\
\nonumber
&&B_{HPZ}(t) = 2 \int_0^t c(s) f(s) g(s) \, ds -  \int_0^t d(s) \big( f(s)^2 + \widetilde{f}(s) g(s)  \big) \, ds \\ \nonumber && = \frac{\delta \pi \Gamma}{24 \beta} \left\{ 12 \int_0^t s \Big(\big(4 f(s) + s \widetilde{f}(s)\big) g(s) + s f(s)^2 \Big) \, ds \right. \\ \nonumber &&\left. - 8 \Gamma \int_0^t s^2 \Big(\big(3 f(s) + s \widetilde{f}(s)\big) g(s) + s f(s)^2 \Big) \, ds + \Gamma^2 \int_0^t s^3 \Big(\big(8 f(s) + 3s \widetilde{f}(s)\big) g(s) + 3 s f(s)^2 \Big)  \, ds \right. \\ \label{B-HPZ} && \left. - \Omega^2 \int_0^t s^3 \Big(\big(8 f(s) + s \widetilde{f}(s)\big) g(s) + s f(s)^2 \Big)  \, ds\right\} \, , \\
\nonumber
&&C_{HPZ}(t) = \int_0^t c(s) f(s)^2 \, ds -  \int_0^t d(s) f(s) \widetilde{f}(s) \, ds \\  \nonumber && = \frac{\delta \pi \Gamma}{24 \beta} \left\{ 12 \int_0^t s f(s) \big(2 f(s) + s \widetilde{f}(s)\big) \, ds -4 \Gamma \int_0^t s^2 f(s) \big(3 f(s) +2 s \widetilde{f}(s)\big) \, ds \right. \\ \label{C-HPZ}  && \left. + \Gamma^2  \int_0^t s^3 f(s) \big(4 f(s) +3s\widetilde{f}(s)\big) \, ds - \Omega^2 \int_0^t s^3 f(s) \big(4 f(s) +s \widetilde{f}(s)\big) \, ds  \right\} \, .
%  %= \frac{8 \gamma^2}{\pi} \big( I_1'(0) - 2 I_1'(t) e^{2\gamma t} + I_1'(0) e^{4\gamma t}\big) \, ,
\end{eqnarray}
As consequence, the change of variables $z=\widetilde{f}(t) \eta + \widetilde{g}(t) y$, $v = f(t) \eta + g(t) y$, leads to
\begin{equation}
\widehat{G_0^{HPZ}}(t, z, v) = e^{- \widetilde{A_{HPZ}}(t) |z|^2 - \widetilde{B_{HPZ}}(t) z \cdot v - \widetilde{C_{HPZ}}(t) |v|^2} \, ,
\label{hatG-HPZ}
\end{equation}
with
\begin{eqnarray}
\label{tildeA-HPZ}
\widetilde{A_{HPZ}}(t) = \frac{f(t)^2 A_{HPZ}(t) -  f(t)g(t) B_{HPZ}(t) + g(t)^2 C_{HPZ}(t)}{\big(f(t)^2 - \widetilde{f}(t) g(t) \big)^2}   \, , &&\\
\label{tildeB-HPZ}
\widetilde{B_{HPZ}}(t) = \frac{- 2 f(t) \widetilde{f}(t)  A_{HPZ}(t) + \big( f(t)^2+ \widetilde{f}(t) g(t)\big) B_{HPZ}(t) -  2f(t)g(t)  C_{HPZ}(t)}{\big(f(t)^2 - \widetilde{f}(t) g(t) \big)^2} \, ,&& \\
\label{tildeC-HPZ}
\widetilde{C_{HPZ}}(t) = \frac{\widetilde{f}(t)^2 A_{HPZ}(t) -  f(t)\widetilde{f}(t) B_{HPZ}(t) + f(t)^2 C_{HPZ}(t)}{\big(f(t)^2 - \widetilde{f}(t) g(t) \big)^2}  \, .&&
\end{eqnarray}

The fundamental solution is then computed as established in the following 

\begin{lemma}[Fundamental solution of the linear Hu-Paz-Zhang master equation]
Assume that $\Omega, \Gamma, \beta > 0$ are sufficiently small in the sense stated above. Then, the time-dependent coefficients $\widetilde{A_{HPZ}}$, $\widetilde{B_{HPZ}}$, and $\widetilde{C_{HPZ}}$ of $\widehat{G_0^{HPZ}}$ are given by
% \begin{eqnarray}
% \label{A}
% A_{UZ}(t) &=& \frac{1}{2 \gamma \pi} I_1'(0) \Big(  4 e^{-2\gamma t} + \left(4\gamma t - 10\right) e^{2\gamma t} + e^{4\gamma t} + 20\gamma t - 4 \gamma^2 t^2+ 5 \Big) \, , \\
% \label{B}
% B_{UZ}(t) &=&  -\frac{2}{\pi} I_1'(0) \Big(  \left(2\gamma t -6\right) e^{2\gamma t} + e^{4\gamma t} + 6\gamma t + 5 \Big) \, ,\\
% \label{C}
% C_{UZ}(t) &=& \frac{2\gamma}{\pi} I_1'(0) (e^{2\gamma t}-1)^2 \, .
% \end{eqnarray}
% As consequence, 
\begin{equation}
\label{tildes}
\widetilde{A_{HPZ}}(t) = \frac{\delta \pi \Gamma}{2 \beta} t^4 \frac{F_A(t)}{\chi(t)} \, ,  \quad \widetilde{B_{HPZ}}(t) = \frac{\delta \pi \Gamma}{2 \beta} t^3 \frac{F_B(t)}{\chi(t)} \, , \quad \widetilde{C_{HPZ}}(t) = \frac{\delta \pi \Gamma}{2 \beta} t^2 \frac{F_C(t)}{\chi(t)} \, ,
\end{equation}
with
\begin{eqnarray*}
F_A(t) &=&  \frac{1}{4} - \frac{\Gamma}{15} t + \left(\frac{\Gamma^2 -3 \Omega^2}{72} \right) t^2 + \frac{\delta \pi \Gamma^2}{12} t^3 \, , \\
F_B(t) &=& 1 - \frac{\Gamma}{3} t + \left( \frac{\Gamma^2 - 3\Omega^2}{12} \right) t^2 + \frac{\delta \pi \Gamma^2}{3} t^3  \, , \\
F_C(t) &=&  1 - \frac{\Gamma}{3} t + \left( \frac{\Gamma^2 - 4\Omega^2}{12} \right) t^2 + \frac{\delta \pi \Gamma^2}{3} t^3 \, , \\
\chi(t) &=& \Big( 1+ \frac{\delta \pi \Gamma^2}{12} t^3 (2-\Gamma t) \Big)^2 \, .
\end{eqnarray*}
Then, the fundamental solution of the linear HPZ master equation (\ref{HPZeq}) is given by 
\begin{equation}
G_{HPZ}(t,x,\xi,z,v) = G_0^{HPZ}\left(t, x - \nu(t) v - \kappa(t) z, \xi - \widetilde{\nu}(t) v - \widetilde{\kappa}(t) z \right) \, ,
\label{G-HPZ}
\end{equation}
where
\begin{eqnarray*}
&&\nu(t) = t - \frac{\Omega^2}{6} t^3 + \frac{\delta \pi \Gamma^2}{24} t^4 -  \frac{\delta \pi \Gamma^3}{120} t^5  \, , \\  &&\kappa(t) = \widetilde{\nu}(t) = 1 - \frac{\Omega^2}{2} t^2 + \frac{\delta \pi \Gamma^2}{6} t^3 - \frac{\delta \pi \Gamma^3}{24} t^4 \, , \\  &&\widetilde{\kappa}(t) = - \Omega^2 t +  \frac{\delta \pi \Gamma^2}{2} t^2 - \frac{\delta \pi \Gamma^3}{6} t^3 \, . 
\end{eqnarray*}
Here,
\begin{equation}
G_0^{HPZ}(t,x,\xi) = d_{HPZ}(t) \, e^{-a_{HPZ}(t)|x|^2 + b_{HPZ}(t) x \cdot \xi - c_{HPZ}(t) |\xi|^2} \, , 
\label{G0-HPZ}
\end{equation}
with
\begin{equation}
d_{HPZ}(t) = \frac{1}{(2\pi)^3} \big( 4 \widetilde{A_{HPZ}}(t) \widetilde{C_{HPZ}}(t) - \widetilde{B_{HPZ}}(t)^2 \big)^{-\frac{3}{2}}
\label{d-HPZ}
\end{equation}
and
\begin{eqnarray}
\label{a1}
a_{HPZ}(t) &=& 4\pi^2 d_{HPZ}(t)^{\frac{2}{3}} \widetilde{C_{HPZ}}(t)  \, ,  \\  \label{b1} b_{HPZ}(t) &=& 4\pi^2 d_{HPZ}(t)^{\frac{2}{3}} \widetilde{B_{HPZ}}(t) \, , \\ \label{c1} c_{HPZ}(t) &=& 4\pi^2 d_{HPZ}(t)^{\frac{2}{3}}  \widetilde{A_{HPZ}}(t) \, .
\label{abc}
\end{eqnarray}
\end{lemma}
\noindent {\bf{Proof.}} 
Formulae (\ref{A-HPZ})-(\ref{C-HPZ}) lead to the following expressions for the coefficients $A_{HPZ}(t)$, $B_{HPZ}(t)$, and $C_{HPZ}(t)$:
\begin{eqnarray}
\label{A}
A_{HPZ}(t) &=& \frac{\delta \pi \Gamma}{2 \beta} t^4 \left( \frac{1}{4}  - \frac{\Gamma}{15} t + \frac{\Gamma^2 - 3 \Omega^2}{72}  t^2 + \frac{\delta \pi \Gamma^2}{24} t^3 \right) \, , \\
\label{B}
B_{HPZ}(t) &=& \frac{\delta \pi \Gamma}{2 \beta} t^3 \left( -1 + \frac{\Gamma}{3} t - \frac{\Gamma^2 - 3 \Omega^2}{12}  t^2 - \frac{\delta \pi \Gamma^2}{6} t^3 \right) \, , 
\\ \label{C}
C_{HPZ}(t) &=& \frac{\delta \pi \Gamma}{2 \beta} t^2 \left( 1 - \frac{\Gamma}{3} t + \frac{\Gamma^2 - 4 \Omega^2}{12}  t^2 + \frac{\delta \pi \Gamma^2}{6} t^3 \right) \, .
\end{eqnarray}
The coefficients with tilde are then directly obtained from (\ref{A})-(\ref{C}). 
Then, by applying the inverse Fourier transform to $ \widehat{G_0^{HPZ}}$ (cf. (\ref{hatG-HPZ})) we find (\ref{G0-HPZ}).
% Clearly $\widetilde{C_{UZ}} > 0$, and therefore $ a_{UZ} > 0$. Also, since  $I_1'(0) =  \frac{1}{2} \ln\left( 1 + \frac{\Gamma^2}{4\gamma^2}\right) > 0$, it suffices to show that $F_1>0$ in order to conclude that $\widetilde{A_{UZ}} > 0$. Indeed, $F_1(0) = 0$ and 
% $$
% F_1\left(\frac{1}{2\gamma}\right) = \frac{1}{4\gamma} \left( \frac{16}{e} + \frac{1}{e^2} -5\right) > 0 
% $$
% for arbitrary $\gamma > 0$. 
% Hence, $F_1(t) > 0$ for all $t>0$. Finally, 
We have just to check that $4 \widetilde{A_{HPZ}}(t) \widetilde{C_{HPZ}}(t) - \widetilde{B_{HPZ}}(t)^2 > 0$ for all $t>0$ in order to give sense to $ d_{HPZ}(t)$. A straightforward computation starting from (\ref{tildeA-HPZ})-(\ref{tildeC-HPZ}) leads us to the relation 
$$
4 \widetilde{A_{HPZ}}(t) \widetilde{C_{HPZ}}(t) - \widetilde{B_{HPZ}}(t)^2 = \frac{1}{\chi(t)} \big( 4 A_{HPZ}(t) C_{HPZ}(t) - B_{HPZ}(t)^2 \big) \, .
$$ 
Now, inserting the expressions (\ref{A})-(\ref{C}) for the coefficients without tilde into the above equation we obtain
$$
4 \widetilde{A_{HPZ}}(t) \widetilde{C_{HPZ}}(t) - \widetilde{B_{HPZ}}(t)^2 = \frac{\delta^2 \pi^2 \Gamma^3}{60 \beta^2 \chi(t)} t^7  > 0 \, ,
$$
which ends the proof.
% with
%$$
% H(t) = (6 \gamma t - 2 \gamma^2 t^2 -5) e^{4\gamma t} + (4 \gamma t - 4 \gamma^2 t^2 +11) e^{2\gamma t} + e^{-2 \gamma t} - 10\gamma t (1 + \gamma t) - 7 \, .
%$$
%Again, $H(0) = 0$ and 
%$$
% H\left(\frac{1}{2\gamma}\right) = -\frac{1}{2} (5e^{2} + 29) + 12 e + \frac{1}{e} > 0 \, ,
%$$
%thus $H(t) > 0$ for all $t>0$, and consequently $d_{UZ}(t)$ is shown to be well defined.
% 

\hfill 
$\square$

The same properties established in Lemma \ref{lema2} for the fundamental solution of the linear UZ equation are now fulfilled by the HPZ propagator. Given that the proof follows exactly the same steps as before, we omit it here for the sake of concision. Besides, the HPZ counterpart of Lemma \ref{lema3} reads as follows.
  
\begin{lemma}[Rates of time growth/decay near $t=0$ of the coefficients of $G_{HPZ}$] Let $a_{HPZ}(t)$, $b_{HPZ}(t)$, $c_{HPZ}(t)$ and $d_{HPZ}(t)$ be given by (\ref{d-HPZ})-(\ref{c1}) and denote 
$$
D_{HPZ}(t) = 4 A_{HPZ}(t) C_{HPZ}(t) - B_{HPZ}(t)^2 \, . 
$$
Then, their behavior near $t=0$ is as follows:   
$$
D_{HPZ}(t) \sim t^6 \, , \quad  a_{HPZ}(t)  \sim t^{-5} \, , \quad b_{HPZ}(t) \sim t^{-4} \, , \quad c_{HPZ}(t) \sim t^{-2} \, , \quad d_{HPZ}(t) \sim t^{-\frac{21}{2}}  \, .
 $$
 \label{lema5}
\end{lemma}
\noindent {\bf{Proof.}} We have
\begin{eqnarray*}
&& A_{HPZ}(0) = A_{HPZ}'(0) = A_{HPZ}''(0) = A_{HPZ}'''(0) = 0, \quad A_{HPZ}^{(iv)}(0) = \frac{3 \delta \pi \Gamma}{\beta} \, , \\
&& B_{HPZ}(0) = B_{HPZ}'(0) = B_{HPZ}''(0) = 0, \quad B_{HPZ}'''(0) = - \frac{3\delta \pi \Gamma}{\beta}  \, , \\
&& C_{HPZ}(0) = C_{HPZ}'(0) = 0, \quad C_{HPZ}''(0) = \frac{\delta \pi \Gamma}{\beta}  \, , 
\end{eqnarray*}
thus
$$
D_{HPZ}(0) =  D_{HPZ}'(0) = \cdots = D_{HPZ}^{(vi)}(0) = 0, \quad D_{HPZ}^{(vii)}(0) = \frac{84\delta^2 \pi^2 \Gamma^3}{\beta^2} \, ,
$$
from which we may conclude that $D_{HPZ}(t)$ goes as $t^6$ in the vicinity of $t=0$.

Also,
$$
\left( \frac{1}{a_{HPZ}}\right)(t) = \frac{4 \widetilde{A_{HPZ}}(t) \widetilde{C_{HPZ}}(t) - \widetilde{B_{HPZ}}(t)^2}{\widetilde{C_{HPZ}}(t)} = \frac{\delta \pi \Gamma^2}{30\beta} \frac{t^5}{F_C(t)} \, ,
$$
that vanishes at $t=0$. As a matter of fact, the first nonvanishing derivative of $\frac{1}{a_{HPZ}}$ at $t=0$ is that of fifth order, $\left( \frac{1}{a_{HPZ}}\right)^{(v)}(0) = \frac{4\delta \pi \Gamma^2}{\beta}$, so that $\frac{1}{a_{HPZ}}$ goes as $t^5$ at $t=0$. 
Analogously,
$$
\left( \frac{1}{b_{HPZ}}\right)(t) =  \frac{4 \widetilde{A_{HPZ}}(t) \widetilde{C_{HPZ}}(t) - \widetilde{B_{HPZ}}(t)^2}{\widetilde{B_{HPZ}}(t)} =  \frac{\delta \pi \Gamma^2}{30\beta} \frac{t^4}{F_B(t)} \, ,
$$
that leads to $\left( \frac{1}{b_{HPZ}}\right)(0) = \left( \frac{1}{b_{HPZ}}\right)'(0) =  \left( \frac{1}{b_{HPZ}}\right)''(0) = \left( \frac{1}{b_{HPZ}}\right)'''(0) = 0$ and 
$ \left( \frac{1}{b_{HPZ}}\right)^{(iv)}(0) = \frac{4\delta \pi \Gamma^2}{5\beta}$, as well as
$$
\left( \frac{1}{c_{HPZ}}\right)(t) = \frac{4 \widetilde{A_{HPZ}}(t) \widetilde{C_{HPZ}}(t) - \widetilde{B_{HPZ}}(t)^2}{\widetilde{A_{HPZ}}(t)} = \frac{\delta \pi \Gamma^2}{30\beta} \frac{t^3}{F_A(t)} \, ,
$$
whose evaluations at $t=0$ give $\left( \frac{1}{c_{HPZ}}\right)(0) =  \left( \frac{1}{c_{HPZ}}\right)'(0) =  0$, $\left( \frac{1}{c_{HPZ}}\right)''(0) = \frac{4\delta \pi \Gamma^2}{5\beta}$.
Finally,
$$
\left( \frac{1}{d_{HPZ}}\right)^{\frac{2}{3}}(t) =  \frac{\delta^2 \pi^4 \Gamma^3 t^7}{15 \beta^2 \chi(t)}   \, ,
$$
whose first nonvanishing derivative is that of seventh order: $\big(d_{HPZ}^{-\frac{2}{3}}\big)^{(vii)}(0) = \frac{336 \delta^2 \pi^4 \Gamma^3}{\beta^2}$.
Now we are done with the proof.

 \hfill 
$\square$
% 

%%%%%%%%%%%%%%%%%%%%%%%%%%%%%%%%%

\section{The nonlinear problems: Local solvability and regularity}

%%%%%%%%%%%%%%%%%%%%%%%%%%%%%%%%%%%%%%%%%%%

This section is devoted to explore the existence of a local-in-time solution to Eq. (\ref{mild}) through a fixed point argument of Banach type. We start by introducing some conventional notation.

Let $T > 0$ and consider the Banach space $X_T = C([0,T]; L^{1,1} \cap L^{1,2})$ endowed with the norm
$$
\|w\|_T  := \sup_{0 \leq t \leq T} \Big\{ \|w(t)\|_{L^{1,1}} + \|w(t)\|_{L^{1,2}} \Big\} \, .
$$
Define also the following closed and bounded subset of $X_T$:
$$
X_T^K = \big\{ w \in X_T : \ w(0,x,\xi) = w_0(x,\xi) \mbox{ a.e.}, \ \|w\|_T \leq K \big\} \, ,
$$
as well as the map $\Gamma: X_T^K \rightarrow X_T$ by means of
$$
\Gamma(w)(t) = G(t)[w_0] -\int_0^t G(t-s)[(H \ast_\xi w)(s)] \, ds \, ,
$$ 
where we denoted $G(t)[f]$ the action of the kinetic Fokker-Planck flux operator $G$ on $f \in L^{1,1}$, namely
$$
G(t)[f] :=  \int_{\Rset^6} G(t, x,\xi,z,v) f(z,v) \, d(z,v) \,.
$$ 
In what follows we generically denote $G, G_0, w, H, V, n$ the already defined magnitudes when the calculations in which they are involved remain the same for both the UZ and the HPZ problems. Otherwise, we use the subscripts UZ or HPZ to identify the equation we are referring to. 

We first show that $\Gamma$ is well--defined. To this aim, some a priori estimates are collected in the following result.

\begin{proposition}
The following estimates are fulfilled:
\begin{itemize}
\item[(i)] For any $1 \leq p,q < \infty$ and $1 \leq m \leq l \leq \infty$ such that $1 + \frac{1}{p} = \frac{1}{r} + \frac{1}{l}$, $1 + \frac{1}{q} = \frac{1}{s} + \frac{1}{m}$, 
\begin{eqnarray*}
\| G_{UZ}(t)[f] \|_{L^{q,p}} &\leq& e^{6\gamma (1 - \frac{1}{m}) t} \big\|G_0^{UZ} \big\|_{L^{s,r}} \|f\|_{L^{m,l}} \, , \\
\| G_{HPZ}(t)[f] \|_{L^{q,p}} &\leq& |\kappa(t)|^{-3(\frac{1}{l} + \frac{1}{m})} |\kappa(t)^2 - \nu(t) \widetilde{\kappa}(t)|^{-3(1 - \frac{1}{l} - \frac{1}{m})} \big\|G_0^{HPZ} \big\|_{L^{s,r}} \|f\|_{L^{m,l}} \, .
\end{eqnarray*}
\item[(ii)] There exists $C>0$ such that the convolution kernel of the nonlinear term defined in (\ref{H}) satisfies 
\begin{eqnarray*}
\|H\|_{L^1(\Rset^3_\xi)} &\leq& C \big(\|n(t)\|_{L^1(\Rset^3)} + \|n(t)\|_{L^2(\Rset^3)}\big) \, , \\
\|\xi H\|_{L^1(\Rset^3_\xi)} &\leq& C \big(\|n(t)\|_{L^1(\Rset^3)} + \|\nabla_x n(t)\|_{L^2(\Rset^3)}\big) \, ,
\end{eqnarray*}
for both the UZ and the HPZ systems. 
%Here, $n(t,x)$ denotes the position density introduced in (\ref{n}).
\end{itemize}
\label{apriori}
\end{proposition}

\noindent {\bf{Proof.}}
The first assertion follows from the change of variables $z + \frac{1}{2\gamma} (1 - e^{-2\gamma t}) v \mapsto z$, then $e^{-2\gamma t} v \mapsto v$ for the UZ case; and $\kappa(t) z + \nu(t) v \mapsto z$, then $\Big( \widetilde{\nu}(t) - \frac{\nu(t) \widetilde{\kappa}(t)}{\kappa(t)}  \Big) v + \widetilde{\kappa}(t) z \mapsto v$ for the HPZ case. The proof concludes after application of Young's inequality for the resulting convolution and Minkowski's inequality for the corresponding norm of $f$.

To prove (ii) we notice that $|H(t, x,\xi)| \leq 16 \big|\check{V}(t,\xi)\big|$, thus we need to control 
$$
\| \check{V}(t)\|_{L^1(\Rset^3)} = \left\| \frac{\check{n}(t, \cdot)}{|\cdot|^2} \right\|_{L^1(\Rset^3)}.
$$
To that aim, we proceed to estimate the $L^1$ norm outside and inside the unit ball $B \subset \Rset^3$, separately. We find 
\begin{eqnarray*}
\left\| \frac{\check{n}(t, \cdot)}{|\cdot|^2} \right\|_{L^1(\Rset^3 \setminus B)} &\leq& C \|\check{n}(t)\|_{L^{2}(\Rset^3)} \leq C  \|n(t)\|_{L^{2}(\Rset^3)} \, , \\
\left\| \frac{\check{n}(t, \cdot)}{|\cdot|^2} \right\|_{L^1(B)} &\leq& C \|\check{n}(t)\|_{L^{\infty}(\Rset^3)} \leq C \, \|n(t)\|_{L^1(\Rset^3)} \, .
\end{eqnarray*}

The estimate for $\xi H$ is analogous, by just considering 
$$
\left\| \frac{\check{n}(t, \cdot)}{|\cdot|} \right\|_{L^1(\Rset^3 \setminus B)}  = \left\| \frac{|\cdot| \check{n}(t, \cdot)}{|\cdot|^2} \right\|_{L^1(\Rset^3 \setminus B)} \leq C  \|\nabla_x n(t)\|_{L^{2}(\Rset^3)} 
$$
by virtue of the differentiability properties of the Fourier transform. This ends the proof.

\hfill
$\square$

As consequence of the a priori estimates established in Proposition \ref{apriori}, the following result holds:

\begin{corollary}
Given $1 \leq p < \infty$, the following $L^{1,p}$ estimates are fulfilled for the action of $G_{UZ/HPZ}$:
\begin{itemize}
\item[(i)] $\|G_{UZ}(t)[w_0]\|_{L^{1,p}} \leq \|w_0\|_{L^{1,p}}$, 
\item[(ii)] $\|G_{UZ}(t-s)[(H \ast_\xi w)(s)]\|_{L^{1,p}} \leq C \|w_{UZ}\|_T \|w_{UZ}(s)\|_{L^{1,p}}$, 
\item[(iii)] $\|G_{HPZ}(t)[w_0]\|_{L^{1,p}} \leq |G(t)|^{-\frac{3}{p'}} |\Lambda(t)|^{\frac{3}{p}} \|w_0\|_{L^{1,p}}$,
\item[(iv)] $\|G_{HPZ}(t-s)[(H \ast_\xi w)(s)]\|_{L^{1,p}} \leq |\kappa(t-s)|^{-\frac{3}{p'}} |\Lambda(t-s)|^{\frac{3}{p}} \|w_{HPZ}\|_T  \|w_{HPZ}(s)\|_{L^{1,p}}$,
\end{itemize}
%\|G_{UZ}(t)[w_0]\|_{L^{1,2}} \leq \|w_0\|_{L^{1,2}} \, , \\
%\|G_{UZ}(t-s)[(H \ast_\xi w)(s)]\|_{L^{1,2}} \leq C T \|w\|_T^2 \, , \\
%\|G_{HPZ}(t)[w_0]\|_{L^{1,1}} \leq |\Lambda(t)|^3 \|w_0\|_{L^{1,1}} \, , \\
%\|G_{HPZ}(t-s)[(H \ast_\xi w)(s)]\|_{L^{1,1}} \leq  C \Big( \int_0^t |\Lambda(t-s)|^3 \, ds \Big) \|w\|_T^2 \, , \\
%\|G_{HPZ}(t)[w_0]\|_{L^{1,2}} \leq  \left| \frac{\Lambda(t)}{G(t)} \right|^{\frac{3}{2}} \|w_0\|_{L^{1,2}} \, , \\
%\|G_{HPZ}(t-
%\end{itemize}
where we denoted $\Lambda(t) := 1 - \frac{\nu(t) \widetilde{\kappa}(t)}{\kappa(t)^2}$. In the particular case $p=1,2$, the right-hand side of (ii) and (iv) can be reformulated as follows:
\begin{itemize}
\item[(ii')] $\|G_{UZ}(t-s)[(H \ast_\xi w)(s)]\|_{L^{1,p}} \leq C T \|w\|_T^2$,
\item[(iv')] $\|G_{HPZ}(t-s)[(H \ast_\xi w)(s)]\|_{L^{1,p}} \leq |\kappa(t-s)|^{-\frac{3}{p'}} |\Lambda(t-s)|^{\frac{3}{p}} T \|w\|_T^2$.
\end{itemize} 
\label{cor}
\end{corollary}

We are now prepared to rigorously set the problem in the context of the Banach fixed point theorem.  

\begin{lemma}
There exist $T, K > 0$ such that $\Gamma: X_T^K \rightarrow X_T$ maps $X_T^K$ onto itself in a contractive way.
\label{Gamma}
\end{lemma}

\noindent {\bf{Proof.}} Take $w \in X_T^K$. On one hand, we have
\begin{eqnarray*}
\|\Gamma_{UZ}[w]\|_{L^{1,1}} &\leq& \|G_{UZ}(t)[w_0]\|_{L^{1,1}}  + \|G_{UZ}(t-s)[(H \ast_\xi w)(s)]\|_{L^{1,1}}  \\ &\leq& \|w_0\|_{L^{1,1}}  + CT K^2 \, , 
\end{eqnarray*}
\begin{eqnarray*}
\|\Gamma_{UZ}[w]\|_{L^{1,2}} &\leq& \|G_{UZ}(t)[w_0]\|_{L^{1,2}} + \|G_{UZ}(t-s)[(H \ast_\xi w)(s)]\|_{L^{1,2}} \\ &\leq& \|w_0\|_{L^{1,2}}  + CTK^2 \, ,
\end{eqnarray*}
and $K = 2 (\|w_0\|_{L^{1,1}} + \|w_0\|_{L^{1,2}})$, $T \leq \frac{1}{4CK}$ can be chosen to find $\|\Gamma_{UZ}[w]\|_T \leq K$.

On the other hand, regarding the HPZ estimates in Corollary \ref{cor} (v) and (vi) it is enough to consider, for instance, a value $T^\ast < \sqrt{\frac{\sqrt{\Omega^2 + \delta \pi \Gamma} - \Omega}{\delta \pi \Gamma}}$ in order to avoid the vertical asymptote of the function $\Lambda(t)$ (reached at $t= \sqrt{2}\sqrt{ \frac{\sqrt{3} \sqrt{2\pi \delta \Gamma^3 + 3 \Omega^4} - 3 \Omega^2}{\pi \delta \Gamma^3}}$) and guarantee that $\Lambda(t) < 2$ for all $0 \leq t \leq T^\ast$ (cf. Figure 1). Indeed, we have
$$
0 \leq 1 - 2 \Omega t^2 - \delta \pi \Gamma t^4 < 1 - 2 \Omega^2 t^2 + \frac{5 \delta \pi \Gamma^2}{6} t^3 - \frac{\delta \pi \Gamma^3}{4} t^4 = \nu(t) \widetilde{\kappa}(t) + \kappa(t)^2 
$$
if and only if $0 \leq t \leq T^\ast$.
Hence,  
\begin{eqnarray*}
\|\Gamma_{HPZ}[w]\|_{L^{1,1}}  &\leq& \|G_{HPZ}(t)[w_0]\|_{L^{1,1}}  + \|G_{HPZ}(t-s)[(H \ast_\xi w)(s)]\|_{L^{1,1}}  \\ &\leq& |\Lambda(t)|^3 \|w_0\|_{L^{1,1}}  +  C \Big( \int_0^t |\Lambda(t-s)|^3 \, ds \Big)K^2  \\ &\leq&  8 \|w_0\|_{L^{1,1}}  +  8C T^\ast K^2 \, .
\end{eqnarray*}
The $L^{1,2}$ norm is estimated analogously, by now choosing $T^{\ast \ast} = \sqrt{3} T^\ast$ so that $\big| \frac{\Lambda(t)}{\kappa(t)} \big| < 4$ for all $0 \leq t \leq T^{\ast \ast}$. Indeed,
$$
0 \leq 3 - 6 \Omega t^2 - 3 \delta \pi \Gamma t^4 < 3 - 6 \Omega^2 t^2 + \frac{13 \delta \pi \Gamma^2}{6} t^3 - \frac{7 \delta \pi \Gamma^3}{4} t^4 = \nu(t) \widetilde{\kappa(t)} - \kappa(t)^2 + 4 \kappa(t)^3 
$$
if and only if $0 \leq t \leq T^{\ast \ast}$. Hence,
\begin{eqnarray*}
\|\Gamma_{HPZ}[w]\|_{L^{1,2}} &\leq& \|G_{HPZ}(t)[w_0]\|_{L^{1,2}} + \|G_{HPZ}(t-s)[(H \ast_\xi w)(s)]\|_{L^{1,2}} \\ &\leq&  \left| \frac{\Lambda(t)}{\kappa(t)} \right|^{\frac{3}{2}} \|w_0\|_{L^{1,2}}  + C \Big( \int_0^t \left| \frac{\Lambda(t-s)}{\kappa(t-s)} \right|^{\frac{3}{2}} \, ds \Big) K^2 \\ &\leq& 8 \|w_0\|_{L^{1,2}}  + 8C T^\ast K^2 \, .
\end{eqnarray*}
Then, $K = 16 (\|w_0\|_{L^{1,1}} + \|w_0\|_{L^{1,2}})$ and $T \leq {\hbox{min}}\{ T^\ast, \frac{1}{32CK} \}$ can be chosen to find $\|\Gamma_{HPZ}[w]\|_T \leq K$.

There only remains to prove the contractivity of $\Gamma: X_T^K \rightarrow X_T^K$, but it follows immediately from 
$$
\| \Gamma(w_1) - \Gamma(w_2) \|_T \leq \frac{1}{2} \|w_1 - w_2\|_T 
$$
for both cases, given any $w_1, w_2 \in X_T^K$. Now we are done with the proof.

\hfill
$\square$

\begin{figure}[h!] \centering
\includegraphics[scale=.28]{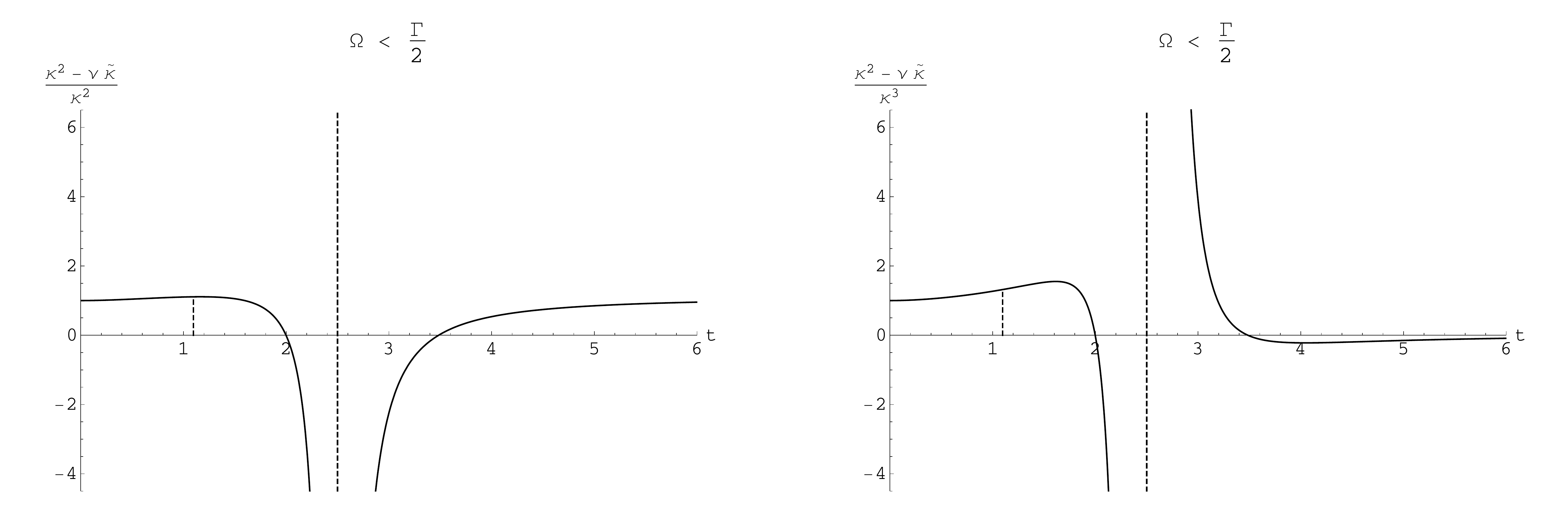} 
\includegraphics[scale=.28]{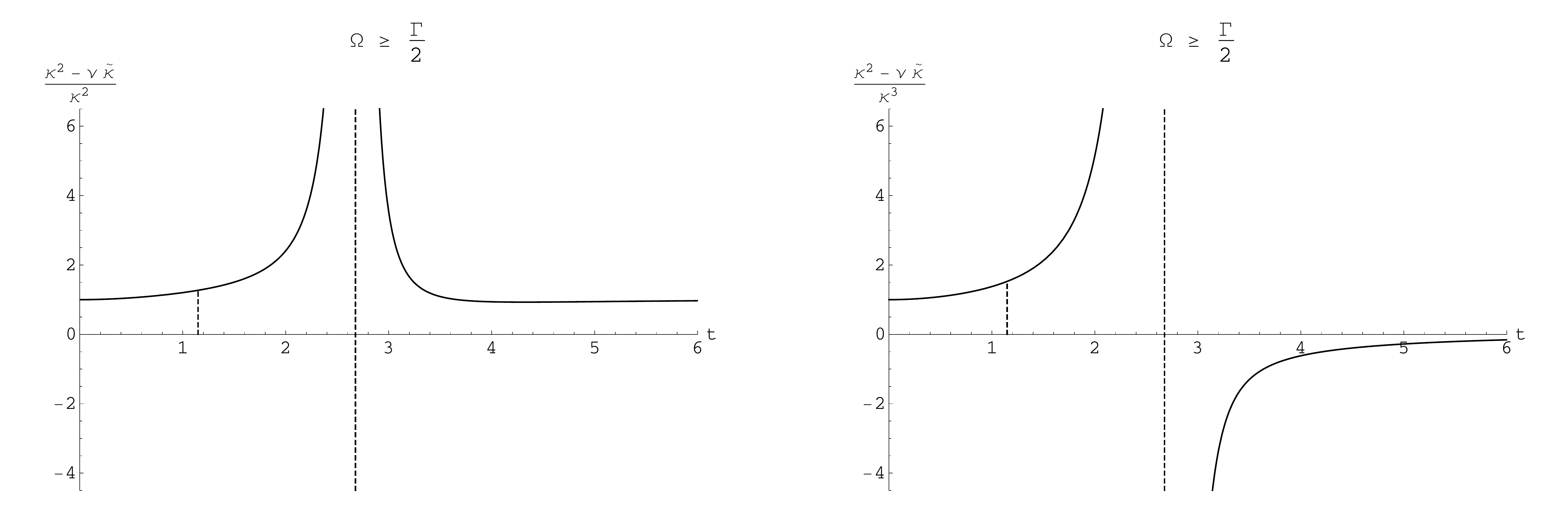} 
\caption{\label{fig1} {Graphical visualization of the functions $\Lambda(t)$ (left-hand side) and 
$\frac{\Lambda(t)}{\kappa(t)}$ (right-hand side) in different situations ($\Omega < \frac{\Gamma}{2}$ and $\Omega \geq \frac{\Gamma}{2}$). 
The short dashed lines correspond to the upper bound chosen in the proof of Lemma \ref{Gamma}: 
$t = \sqrt{\frac{\sqrt{\Omega^2 + \delta \pi \Gamma} - \Omega}{\delta \pi \Gamma}}$.}}
\end{figure}

As consequence, it can be claimed that $\Gamma$ has a unique fixed point $w \in X_T^K$, which is equivalent to ensure
the existence of a unique mild solution (defined on $[0,T]$ for sufficiently small $T>0$ depending on $w_0$) of the 
initial value problem associated with the UZ and the HPZ master equations. By virtue of Pazy's theory \cite{Pazy}, 
there is a maximum time of existence $T_{max}$ which is either $T_{max} = \infty$, or $T_{max} < \infty$ such that 
$\lim_{T \to T_{max}} \{ \|w\|_T \} = \infty$. As a matter of fact, the next section is devoted to prove that the second 
possibility cannot occur, thus $w(t,x,\xi)$ is defined globally in time. We close this section with some regularity properties 
regarding the Wigner function, the position density and the potential.

\begin{proposition}[Regularity]
Let $0<T<T_{max}$. Then, the following properties are satisfied:
%and let also $w_{UZ}$ and $w_{HPZ}$ be the mild solutions of the UZ and the HPZ initial value problems, respectively. Then
 \begin{itemize}
  \item[(i)] $w \in C((0,T); L^{\infty, \infty}) \cap C((0,T); L^{1,p})$, for all $1 \leq p < \infty$. 
  \item[(ii)] $\nabla_x V(t) \in L^2(\Rset^3) \cap L^\infty(\Rset^3)$.
 % \item[(iii)] $\nabla_x n(t) \in L^2(\Rset^3)$.
%  \item[(iii)] $\xi w(t) \in L^{1,1}$ for all $0<t<T$. As a matter of fact, the inequalities
%  $$
%  \|\xi w_{UZ}(t)\|_{L^{1,1}} \leq C_1 e^{C_2 t} \, , \quad  \|\xi w_{HPZ}(t)\|_{L^{1,1}} \leq 
%  $$
%  hold for all $0<t<T$, where $C_1$ and $C_2$ are positive constants depending on $\|w\|_T$.
 \end{itemize}
\label{reg}
\end{proposition}

\noindent {\bf{Proof.}} We prove (i) in several steps. The first one consists in using Proposition \ref{apriori} (i) and Lemma \ref{lema2} (ii) with $q=p$ to estimate
\begin{eqnarray}
\nonumber
||w_{UZ}(t)||_{L^{p_1, p_1}} &\leq& ||G_{UZ}(t)[w_0]||_{L^{p_1, p_1} }+ \int_0^t||G_{UZ}(t-s)[(H\ast_\xi w_{UZ})(s)]||_{L^{p_1, p_1}} \, ds \\
&\leq& C ||w_0||_{L^{1,1}} \, d_{UZ}(t)^{ \frac{1}{p'_1} } + C ||w_{UZ}||_T^2 \int_0^t d_{UZ}(s)^{ \frac{1}{p'_1} } \, ds \, ,
\label{step1}
\end{eqnarray}
which entails $w\in C((0,T_{max}); L^{p_1}(\Rset^6))$ for all $1\leq p_1 < \frac{9}{8}$. Indeed, since $d_{UZ}\stackrel {t=0} {\sim}t^{-9}$ by virtue of Lemma \ref{lema3}, then $d_{UZ}^{\frac{1}{p'_1}}\stackrel {t=0} {\sim}t^{\frac{9}{p_1}-9}$ remains integrable for the aforesaid range of values for $p_1$.
%\[\begin{aligned}&d_{UZ}\stackrel {t=0} {\sim}t^{-9}\Rightarrow d_{UZ}^{\frac{1}{p'}}\stackrel {t=0} {\sim}t^{-\frac{9}{p'}}=t^{\frac{9}{p}-9}\Rightarrow\\
%&\frac{9}{p}-9>-1 \Rightarrow  9-9p>-p \Rightarrow p<\frac{9}{8}.
%  \end{aligned}
In a second step we start from an arbitrary time $\varepsilon>0$ and rewrite $w_{UZ}$ as follows:

\[
 w_{UZ}(t) = G_{UZ}(t-\varepsilon)[w(\epsilon)] - \int_{\varepsilon}^tG_{UZ}(t-s)[(H\ast_\xi w_{UZ})(s)] \, ds \, ,
\]
\noindent
so that an estimate of $||w_{UZ}(t)||_{L^{p_2, p_2}}$ like that in (\ref{step1}) can be performed for all $p_2 < \frac{9}{7}$ and $t > \varepsilon$, as follows from the identifications $q=p=p_2$, $l=m=p_1$ and $r=s < \frac{9}{8}$ in Proposition \ref{apriori} (i). Proceeding analogously,  after the $j$-th step we obtain an estimate for $||w_{UZ}(t)||_{L^{p_j, p_j}}$ with $p_j < \frac{9}{9-j}$ and $t > (j-1)\varepsilon$. After nine steps we find 
$$
w_{UZ}(t) = G_{UZ}(t-8\varepsilon)[w_{UZ}(8\varepsilon)] - \int_{8\varepsilon}^t G_{UZ}(t-s)[(H\ast_\xi w_{UZ})(s)] \, ds \, ,
$$
which leads us to an estimate of $||w_{UZ}(t)||_{L^{p_9, p_9}}$ for all $p_9<\infty$ and $t>8\epsilon$. Finally, in a tenth step we get a uniform bound for $w_{UZ}$ as follows:
\begin{eqnarray*}
||w_{UZ}(t)||_{L^{\infty, \infty}} &\leq& C d_{UZ}(t-9\varepsilon)^{\frac{1}{p'_{9}}} ||w_{UZ}(9\varepsilon)||_{L^{r,r}} \\
&& + \, C ||w_{UZ}||_T \int_{9\varepsilon}^t d_{UZ}(t-s)^{\frac{1}{p'_{9}}}||w_{UZ}(s)||_{L^{r,r}} \, ds \, ,
\end{eqnarray*}
for all $t>9\varepsilon$. From the arbitrariness of $\varepsilon$ we conclude (i) for the UZ system. 

The HPZ case is analogous. In a first step we obtain
\begin{eqnarray*}
  ||w_{HPZ}(t)||_{L^{p_1, p_1}} &\leq& ||G_{HPZ}(t)[w_0]||_{L^{p,p}} + \int_0^t ||G_{HPZ}(t-s)[(H\ast_\xi w_{HPZ})(s)]||_{L^{p,p}} \, ds \\
  &\leq& C ||w_0||_{L^{1,1}} \, d_{HPZ}(t)^{\frac{1}{p'_1}} + C ||w_{HPZ}||_T^2 \int_0^t d_{HPZ}(s)^{\frac{1}{p'_1}} \, ds \, .
 \end{eqnarray*}
Since $d_{HPZ} \stackrel {t=0} {\sim}t^{-\frac{21}{2}}$ by virtue of Lemma \ref{lema5}, then $d_{HPZ}^{\frac{1}{p'_1}}\stackrel {t=0} {\sim}t^{\frac{21}{2 p_1} - \frac{21}{2}}$ remains integrable for all $1 \leq p_1 < \frac{21}{19}$. The second step leads to
\begin{eqnarray*}
||w_{HPZ}(t)||_{L^{p_2, p_2}} &\leq& C d_{HPZ}(t-\varepsilon)^{\frac{1}{p_1'}} ||w_{HPZ}(\varepsilon)||_{L^{p_1, p_1}} \\ &&+ \, C ||w_{HPZ}||_T \int_{\varepsilon}^td_{HPZ}(t-s)^{\frac{1}{p_1'}} ||w_{HPZ}(s)||_{L^{p_1, p_1}} \, ds \, ,
\end{eqnarray*}
which allows for a range of integrability $1 \leq q \leq \frac{21}{17}$. The procedure concludes after twelve steps, when we find the following uniform estimate for $w_{HPZ}$:
\begin{eqnarray*}
 ||w_{HPZ}(t)||_{L^{\infty, \infty}} &\leq& C d_{HPZ}(t -11 \varepsilon)^{\frac{1}{p_{11}'}} ||w_{HPZ}(11\varepsilon)||_{L^{p_{11}, p_{11}}} \\ && 
  + \, C ||w_{HPZ}||_T\int_{11\varepsilon}^td_{HPZ}(t-s)^{\frac{1}{p_{11}'}}||w(s)||_{L^{p_{11}, p_{11}}}ds,
\end{eqnarray*}
for any $t >11\varepsilon$. Again, from the arbitrariness of $\varepsilon$ we deduce (i) in the $HPZ$ case.

 On the other hand, when fixing the exponent of $\xi$-integrability as $q=1$ we find
 \begin{eqnarray*}
%\nonumber
||w_{UZ}(t)||_{L^{1, p}} &\leq& ||G_{UZ}(t)[w_0]||_{L^{1, p}} + \int_0^t||G_{UZ}(t-s)[(H\ast_\xi w_{UZ})(s)]||_{L^{1, p}} \, ds \\
&\leq& ||w_0||_{L^{1,p}} +  ||w_{UZ}||_T \int_0^t  ||w_{UZ}(s)||_{L^{1,p}} \, ds \, ,
\end{eqnarray*}
which leads to $||w_{UZ}(t)||_{L^{1, p}} \leq ||w_0||_{L^{1,p}} \, e^{||w_{UZ}||_T \, t}$ via Gronwall's inequality.
In an analogous way, 
 \begin{eqnarray*}
&&||w_{HPZ}(t)||_{L^{1, p}} \leq ||G_{HPZ}(t)[w_0]||_{L^{1, p}} + \int_0^t||G_{HPZ}(t-s)[(H\ast_\xi w_{HPZ})(s)]||_{L^{1, p}} \, ds \\
&&\leq |\kappa(t)|^{-\frac{3}{p'}} |\Lambda(t)|^{\frac{3}{p}} ||w_0||_{L^{1,p}} +  ||w_{HPZ}||_T \int_0^t  |\kappa(t-s)|^{-\frac{3}{p'}} |\Lambda(t-s)|^{\frac{3}{p}} ||w_{HPZ}(s)||_{L^{1,p}} \, ds \, ,
\end{eqnarray*}
hence
$$
||w_{HPZ}(t)||_{L^{1, p}} \leq C ||w_0||_{L^{1,p}} \, e^{ C ||w_{HPZ}||_T \, T} \, ,
$$
for sufficiently small $T$, as follows from a similar argument to that held in the proof of Lemma \ref{Gamma}.
 
 (ii) The $L^2$ and $L^\infty$ bounds are derived from the general theory of singular integrals of convolution type (see \cite{Stein}). Indeed, we have
 \begin{eqnarray*}
 \| \nabla_x V(t)\|_{L^2(\Rset^3)} &\leq& C \|n(t)\|_{L^{\frac{6}{5}}(\Rset^3)} \leq C \|n(t)\|_{L^2(\Rset^3)}^{\frac{5}{11}} \leq C \|w\|_T^{\frac{5}{11}} \, , \\
  \| \nabla_x V(t)\|_{L^\infty(\Rset^3)} &\leq& C \|n(t)\|_{L^r(\Rset^3)}^{\frac{2}{3}(\frac{r}{r-1})} \leq C \|w(t)\|_{L^{1,r}}^{\frac{2}{3}(\frac{r}{r-1})} \, , \quad 3 < r \leq \infty \, ,
 \end{eqnarray*}
 according to (i).
 
% (iii) Taking $x$-derivatives in (\ref{mild}) and then the $L^{1,2}$ norm we find
% \begin{eqnarray*}
% \|\nabla_x w(t)\|_{L^{1,2}} &\leq&  \|\nabla_x G(t)[w_0]\|_{L^{1,2}} + \int_0^t  \|\nabla_x G(t-s)[H\ast w](s)\|_{L^{1,2}} \, ds
% \end{eqnarray*}
% 
% 
% COMPLETAR
% 
% The proof concludes by noticing that $\|\nabla_x n(t)\|_{L^2} \leq \|\nabla_x w(t)\|_{L^{1,2}}$.
%%Indeed, since indexes $q,p,r$ are linked by Young's relations $\frac{1}{q}=\frac{1}{p} +\frac{1}{r} -1$, 
%%we get in our case that $$\frac{1}{q}>\frac{8}{9}-1 +\frac{1}{r} = -\frac{1}{9} +\frac{8}{9}=\frac{7}{9}\Rightarrow q<\frac{9}{7}$$.

%(iii) We now multiply (\ref{mild}) by $\xi$ and take $L^1$ norms to find
%\begin{eqnarray*}
%\|\xi w_{UZ}(t)\|_{L^{1,1}} &\leq&  \|\xi G_{UZ}(t)[w_0]\|_{L^{1,1}} + \int_0^t  \|\xi G_{UZ}(t-s)[H\ast w_{UZ}](s)\|_{L^{1,1}} \, ds
%\end{eqnarray*}
%
%COMPLETAR

\hfill
$\square$

%%%%%%%%%%%%%%%%%%%%%%%%%%%%%%%%%

\section{The nonlinear problems: Global solvability}

%%%%%%%%%%%%%%%%%%%%%%%%%%%%%%%%%%%%%%%%%%%

Our global existence proof for both the UZ and the HPZ initial value problems relies on the following Lieb-Thirring type estimate for the particle density $n(t,x)$ in terms of the kinetic energy
$$
E[w](t) = \frac{1}{2} \int_{\Rset^6} |\xi|^2 w(t,x,\xi) \, d(\xi,x) \, , 
$$
given that the density matrix operator (represented by $\rho(t,x,y)$, cf. (\ref{w})) remains nonnegative along the evolution \cite{A, LP}:
\begin{equation}
\| n(t) \|_{L^p(\Rset^3)} \leq C(p) \, Q^{r}  E[w](t)^{1-r} \, , \\
%&&\Big\| \int_{\Rset^3} \xi w(t,x,\xi) \, d\xi \Big\|_{L^r(\Rset^3)} \leq C(r) Q^{\theta_2} \left( \int_{\Rset^3} |\xi|^2 w(t,x,\xi) \, d\xi \right)^{1-\theta_2} \, ,
\label{LT}
\end{equation}
with $1 \leq p \leq 3$ and $r = \frac{3-p}{2p}$, where we denoted 
$$
Q[f](t) = \int_{\Rset^6} f(t,x,\xi) \, d(\xi, x) \, .
$$ 
If $f=w$, then $Q[w] \equiv Q$ is a conserved quantity that describes the total mass of the system. 
%and 
%$$
%E[w](t) = \frac{1}{2} \int_{\Rset^6} |\xi|^2 w(t,x,\xi) \, d(\xi,x)
%$$ 
%the kinetic energy. 
We also note that the first order velocity moment 
$$
j[f](t,x) = \int_{\Rset^3} \xi f(t,x,\xi) \, d\xi 
$$
stands for the current density of the system when $f=w$, which is linked with the position density through the standard continuity equation 
\begin{equation}
\partial_t n + \nabla_x \cdot j[w] = 0 \, , 
\label{ce}
\end{equation}
as simply deduced by integrating Eqs. (\ref{UZ}) and (\ref{HPZeq}) against $\xi$.
%In what follows we make use of the following Gaussian regularization of the Wigner function, typically known under the name of Husimi transform, in order to show that the kinetic energy is bounded from below:
%$$
%w_H(t,x,\xi) = w(t,x,\xi) \ast_{x,\xi} \pi^{-3} e^{-|x|^2 - |\xi|^2} \, ,
%$$ 
%whose main feature is to be pointwise nonnegative on $\Rset^3_x \times \Rset^3_\xi$. Also, it is a simple matter to observe that the Husimi kinetic energy is connected to that of $w$ through the formula
%\begin{equation}
%E[w_H](t) = E[w](t) + \frac{3Q}{2} \, .
%\label{ehus}
%\end{equation}
%%\quad 1 \leq r \leq \frac{3}{2} \, ,
%%\quad \theta_2 = \frac{3-2r}{2r} \, .

The following result provides us with the required tools to control the $\| \cdot\|_{L^{1,1}} + \| \cdot \|_{L^{1,2}}$ norm of the solutions for all times. 

\begin{proposition}[Properties of the kinetic energy] The following assertions hold true:

\begin{itemize}
\item[(i)] $\int_{\Rset^3} \nabla_x V(t,x) \cdot j[w](t,x) \, dx = - \frac{1}{2} \partial_t \|\nabla_x V(t)\|_{L^2(\Rset^3)}^2$.
\item[(ii)] $\int_{\Rset^3} n(t,x) (x \cdot \nabla_x V(t,x)) \, dx = \frac{1}{2} \|\nabla_x V(t)\|_{L^2(\Rset^3)}^2$.
\item[(iii)] The kinetic energy associated with the action of the fundamental solution on a  function $f \in L^{1,1}$ is given by
\begin{eqnarray*}
E[G_{UZ}(t)[f]] &=& c_{UZ}(t)^{-1} \big( C_1 + C_2 \, d_{UZ}(t)^{-\frac{2}{3}} b_{UZ}(t)^2 \big) Q[f](t)  + e^{-4\gamma t} E[f](t) \, , \\
E[G_{HPZ}(t)[f]] &=& c_{HPZ}(t)^{-1} \big( C_1 + C_2 \, d_{HPZ}(t)^{-\frac{2}{3}} b_{HPZ}(t)^2 \big) Q[f](t)  + {\widetilde{\nu}}(t)^2 E[f](t) \\ && + \, \frac{1}{2} {\widetilde{\kappa}}(t)^2 \int_{\Rset^6} |z|^2 f(t,z,v) \, d(v,z) + {\widetilde{\nu}}(t) {\widetilde{\kappa}}(t)   \int_{\Rset^6}  z \cdot  j[f](t,z)  \, dz  \, ,
\end{eqnarray*}
for some $C_1, C_2 > 0$. In particular, it follows that
\begin{eqnarray*}
E[G_{UZ}(t)[w_0]] &\leq& E[w_0] +  KQ t^2 \, , \quad K>0 \, , \\  E[G_{UZ}(t-s)[(H \ast w)(s)]] &=& - \frac{1}{2}  e^{-4\gamma (t-s)} \partial_s \|\nabla_x V(s)\|_{L^2(\Rset^3)}^2 \, , \\ E[G_{HPZ}(t)[w_0]] &\leq& E[w_0] + KQ t(1+t) \, , \quad K>0 \, , \\ E[G_{HPZ}(t-s)[(H \ast w)(s)]] &=& - \frac{1}{2} {\widetilde{\nu}}(t-s)^2 \partial_s \|\nabla_x V(s)\|_{L^2(\Rset^3)}^2 \\ &&+  \frac{1}{2}  {\widetilde{\nu}}(t-s) {\widetilde{\kappa}}(t-s) \|\nabla_x V(s)\|_{L^2(\Rset^3)}^2 \, .
\end{eqnarray*}
\item[(iv)] $|E[w](t)| < \infty$ for all $0 < t < T$.
\end{itemize}
\label{prop-global}
\end{proposition}

\noindent {\bf{Proof.}}
To prove (i) we just notice that
$$
\int_{\Rset^3} \nabla_x V(t,x) \cdot j[w](t,x) \, dx = \int_{\Rset^3} V(t,x) \partial_t n(t,x) \, dx  = - \frac{1}{2}  \partial_t \|\nabla_x V(t) \|_{L^2(\Rset^3)} \, ,
$$
by virtue of the continuity equation (\ref{ce}) and the Poisson relation $\Delta_x V = n$.

Assertion (ii) shows up after integrating by parts repeatedly:
\begin{eqnarray*}
 \int_x n(t,x) \big( x \cdot \nabla_x V(t,x) \big) \, dx 
%= \int_x \Delta_x V(t,x) \big( x \cdot \nabla_x V(t,x) \big) \,  dx \\ && 
&=& - \int_x \nabla_x V(t,x) \cdot \nabla_x \big( x \cdot \nabla_x V(t,x) \big) \,  dx \\ &=&  - \| \nabla_x V(t)\|_{L^2}^2 - \alpha(t) \, ,
\end{eqnarray*}
with
\begin{eqnarray*}
\alpha(t) &=& \sum_{i,j =1}^3 \int_x \partial_{x_j} V(t,x) \, x_i \,  \partial^2_{x_j x_i}V(t,x)  \, dx \\ &=& - \sum_{i,j =1}^3 \int_x \partial_{x_i} \big( x_i \partial_{x_j} V(t,x) \big)\partial_{x_j}V(t,x)  \, dx = - 3 \| \nabla_x V(t)\|_{L^2}^2 - \alpha(t) \, , 
\end{eqnarray*}
hence $\alpha(t) =  - \frac{3}{2} \| \nabla_x V(t)\|_{L^2}^2$. 

We now prove (iii). A straightforward computation leads to 
\[\begin{aligned}
&E[G_{UZ}(t)[f]]\\
&=\frac{1}{2} \int_{\Rset^6}  \left( \int_{\Rset^6} |\xi + e^{-2 \gamma t}v|^2 G_0^{UZ}\left(t,x - z -
\frac{1}{2\gamma}(1 - e^{-2\gamma t})v,\xi\right) \, d(\xi,x) \right) f(t,z,v) \, d(v,z) \\ 
&= E[G_0^{UZ}](t) Q[f](t) + e^{-4\gamma t}E[f](t) + e^{-2\gamma t}  \int_{\Rset^3} j[G_0^{UZ}](t,x) \, dx \cdot \int_{\Rset^3} j[f](t,z) \, dz \,.
\end{aligned}\]
We have
\begin{eqnarray*}
E[G_0^{UZ}](t) &=& \frac{1}{2} d_{UZ}(t) \int_{\Rset^3} e^{-a_{UZ}(t)|x|^2} \left( \int_{\Rset^3} |\xi|^2 e^{b_{UZ}(t)x \cdot 
\xi - c_{UZ}(t) |\xi|^2} \, d\xi \right) dx \\ &=& \frac{1}{2} d_{UZ}(t) c_{UZ}(t)^{-\frac{5}{2}} \int_x e^{- \frac{d_{UZ}(t)^{\frac{2}{3}} }{4 c_{UZ}(t)} |x|^2} \left( \int_{\Rset^3} \Big|\xi + \frac{b_{UZ}(t)x}{2 \sqrt{c_{UZ}(t)}} \Big|^2 e^{- |\xi|^2 } \, d\xi \right) dx \\ &=& c_{UZ}(t)^{-1} \big( C_1 + C_2 \, d_{UZ}(t)^{-\frac{2}{3}} b_{UZ}(t)^2 \big) 
\end{eqnarray*}
for some $C_1, C_2 > 0$, and
\begin{eqnarray*}
 && \int_{\Rset^3} j[G_0^{UZ}](t,x) \, dx = d_{UZ}(t) \int_{\Rset^3} e^{-a_{UZ}(t)|x|^2} \left( \int_{\Rset^3} \xi \, e^{b_{UZ}(t)x \cdot \xi - c_{UZ}(t) |\xi|^2} \, d\xi \right) dx \\ && = d_{UZ}(t) c_{UZ}(t)^{-2} \int_{\Rset^3} e^{- \frac{d_{UZ}(t)^{\frac{2}{3}} }{4 c_{UZ}(t)} |x|^2} \left( \int_{\Rset^3} \Big(\xi +  \frac{b_{UZ}(t)x}{2 \sqrt{c_{UZ}(t)}}\Big) e^{- |\xi|^2 } \, d\xi \right) dx = 0 \, .
 \end{eqnarray*}
Thus, assertion (iii) follows in the UZ case. For the HPZ system we proceed analogously to deduce
\[\begin{aligned}
&E[G_{HPZ}(t)[f]]\\ 
&= \frac{1}{2} \int_{\Rset^6}  \left( \int_{\Rset^6} |\xi +  {\widetilde{\nu}}(t)v + {\widetilde{\kappa}}(t)z|^2 \, 
G_0^{HPZ}(t, x-\nu(t)v-\kappa(t)z, \xi) \, d(\xi,x) \right) f(t,z,v)  \, d(v,z)\\ 
&= E[G_0^{HPZ}](t) Q[f](t) + {\widetilde{\nu}}(t)^2 E[f](t) + \, \frac{1}{2} {\widetilde{\kappa}}(t)^2 \int_{\Rset^6} |z|^2 f(t,z,v) \, d(v,z)\\  
&+ \,  \int_{\Rset^3}  j[G_0^{HPZ}](t,x) \, dx  \cdot \left( {\widetilde{\nu}}(t) \int_{\Rset^3} j[f](t,z) \, dz 
+ {\widetilde{\kappa}}(t)  \int_{\Rset^6}  z f(t,z,v)  \, d(v,z)\right)
%\\ && + \, {\widetilde{G}}(t)   \int_{\Rset^3}  j[G_0^{HPZ}](t,x) \, dx \cdot \int_{\Rset^6}  z f(t,z,v)  \, d(v,z) 
\\ 
&+ \, {\widetilde{\nu}}(t) {\widetilde{\kappa}}(t)   \int_{\Rset^3}  z \cdot j[f](t,z)  \, dz \, ,
\end{aligned}\]
with
$$
E[G_0^{HPZ}](t) = c_{HPZ}(t)^{-1} \big( C_1 + C_2 \, d_{HPZ}(t)^{-\frac{2}{3}} b_{HPZ}(t)^2 \big) \, , \quad \int_{\Rset^3} j[G_0^{HPZ}](t, x) \, dx = 0 \, ,
$$
exactly as before. The estimates for $E[G_{UZ}(t)[w_0]]$ and $E[G_{HPZ}(t)[w_0]]$ follow from the decay rates given in Lemma \ref{lema3} and \ref{lema5}. Finally, choosing $f = H \ast w$ and taking into account that 
$$
\int_{\Rset^3} H(t,x,\xi) \, d\xi = \int_{\Rset^3} |\xi|^2 H(t,x,\xi) \, d\xi = 0 \, , \quad \int_{\Rset^3} \xi H(t,x,\xi) \, d\xi = \nabla_x V(t,x) \, ,
$$
we find 
\begin{eqnarray*}
E[G_{UZ}(t-s)[(H \ast w)(s)]] &=& e^{-4\gamma (t-s)} E[H\ast w](s) \\ &=& \frac{1}{2} e^{-4\gamma (t-s)} \int_{\Rset^6} \left( \int_{\Rset^3} |\xi + \xi'|^2 H(s,x,\xi) \, d\xi \right) w(s,x,\xi') \, d(\xi',x) \\ &=& - \frac{1}{2} e^{-4\gamma (t-s)} \partial_s \|\nabla_x V(s)\|_{L^2(\Rset^3)}^2 
\end{eqnarray*}
by virtue of (i), as well as
\begin{eqnarray*}
&&E[G_{HPZ}(t-s)[(H \ast w)(s)]] \\&&= {\widetilde{\nu}}(t-s) \left( {\widetilde{\nu}}(t-s) E[H \ast w](s) +  {\widetilde{\kappa}}(t-s)  \int_{\Rset^3}  x\cdot j[H \ast w](s,x) \, dx \right)
%\\&=& \frac{1}{2} {\widetilde{F}}(t-s)^2 \int_x \int_{\xi} \int_{\xi'} |\xi|^2 H(s,x,\xi - \xi') w(s,x,\xi') \, d\xi' \, d\xi \, dx \\&& + \, \frac{1}{2} {\widetilde{G}}(t-s)^2 \int_x \int_{\xi} \int_{\xi'} |x|^2 H(s,x,\xi - \xi') w(s,x,\xi') \, d\xi' \, d\xi \, dx \\ && + {\widetilde{F}}(t-s) {\widetilde{G}}(t-s)  \int_x \int_{\xi}  \int_{\xi'} \xi \cdot x  \, H(s,x,\xi - \xi') w(s,x,\xi') \, d\xi' \, d\xi \, dx 
\\ &&=
- \frac{1}{2} {\widetilde{\nu}}(t-s) \left( {\widetilde{\nu}}(t-s) \, \partial_s \|\nabla_x V(s)\|_{L^2(\Rset^3)}^2 - {\widetilde{\kappa}}(t-s) \, \|\nabla_x V(s)\|_{L^2(\Rset^3)}^2 \right) .
\end{eqnarray*}

Finally, to prove (iv) we first integrate Eq. (\ref{mild}) against $|\xi|^2$ to obtain 
\begin{eqnarray*}
E[w](t) &=& E[G(t)[w_0]] - \int_0^t E[G(t-s)[(H\ast_{\xi} w)(s)]] \, ds \, . 
\end{eqnarray*}
This leads to
\begin{eqnarray*}
|E[w_{UZ}](t)| &\leq& |E[w_0]| +  KQ t^2 + \frac{1}{2} \Big| \int_0^t e^{-4\gamma (t-s)} \partial_s \|\nabla_x V(s)\|_{L^2(\Rset^3)}^2 \, ds \Big| \, , \\
|E[w_{HPZ}](t)| &\leq& |E[w_0]| +  KQ t(1+t) + \frac{1}{2}  \Big| \int_0^t  {\widetilde{\nu}}(t-s)^2 \, \partial_s \|\nabla_x V(s)\|_{L^2(\Rset^3)}^2 \, ds \Big| \\ &&+ \, \frac{1}{2}  \int_0^t  |{\widetilde{\nu}}(t-s)  {\widetilde{\kappa}}(t-s)| \|\nabla_x V(s)\|_{L^2(\Rset^3)}^2 \, ds \, ,
\end{eqnarray*}
by virtue of (ii). 
Now it is enough to notice that 
\begin{eqnarray*}
%&& \int_0^t e^{-4\gamma (t-s)} \int_{\Rset^3} \nabla_x V(s,x) \cdot j(s,x) \, dx \, ds \\ &&= \int_0^t e^{-4\gamma (t-s)} \int_{\Rset^3} V(s,x) \cdot \partial_s n(s,x) \, dx \, ds \\ && 
&&\int_0^t e^{-4\gamma (t-s)} \, \partial_s \|\nabla_x V(s) \|_{L^2(\Rset^3)}^2 \, ds \\ &&= 
 \|\nabla_x V(t) \|_{L^2(\Rset^3)}^2 - e^{-4\gamma t}  \|\nabla_x V(0) \|_{L^2(\Rset^3)}^2 - 4\gamma \int_0^t e^{-4\gamma (t-s)}  \|\nabla_x V(s) \|_{L^2(\Rset^3)}^2 \, ds  \, ,
% \\ && \leq  \frac{1}{2} \|\nabla_x V(t) \|_{L^2(\Rset^3)}^2 \, , 
%\int_0^t  {\widetilde{F}}(t-s)  {\widetilde{G}}(t-s) \int_{\Rset^3} x \, n(s,x) \cdot \nabla_x V(s,x) \, dx \, ds &=&
\end{eqnarray*}
and
\begin{eqnarray*}
&& \int_0^t  {\widetilde{\nu}}(t-s)^2 \, \partial_s \|\nabla_x V(s)\|_{L^2(\Rset^3)}^2 \, ds \\ &&=  \|\nabla_x V(t) \|_{L^2(\Rset^3)}^2 -  {\widetilde{\nu}}(t)^2  \|\nabla_x V(0) \|_{L^2(\Rset^3)}^2 + 2 \int_0^t  {\widetilde{\nu}}(t-s) {\widetilde{\nu}}'(t-s) \|\nabla_x V(s) \|_{L^2(\Rset^3)}^2 \, ds \, . 
%\\ && \leq \frac{1}{2}  \|\nabla_x V(t) \|_{L^2(\Rset^3)} + \int_0^t  {\widetilde{F}}(t-s) {\widetilde{F}}'(t-s) \|\nabla_x V(s) \|_{L^2(\Rset^3)}^2 \, ds \, .
\end{eqnarray*}
Consequently, Proposition \ref{reg} applies to conclude $|E[w_{UZ}](t)|, |E[w_{HPZ}](t)| \leq C(w_0, Q, T)$.
%\begin{eqnarray*}
%|E[w_{UZ}](t)| &\leq& \Big( |E[w_0]| + \frac{1}{2}  \|\nabla_x V(0) \|_{L^2(\Rset^3)}^2 + \frac{1}{2} \|w_{UZ}\|^{\frac{5}{11}}_T  + 2 \gamma \|w_{UZ}\|^{\frac{5}{11}}_T T \Big) +  KQ T^2 \\ &\leq& C(w_0, Q, T) \, , 
%\end{eqnarray*}
%%\int_0^t  {\widetilde{F}}(t-s)  {\widetilde{G}}(t-s) \int_{\Rset^3} x \, n(s,x) \cdot \nabla_x V(s,x) \, dx \, ds &=&
% \begin{eqnarray*}
%&&|E[w_{HPZ}](t)| \leq \\&& \Big( |E[w_0]| +  \frac{1}{2}  \|w_{HPZ} \|^{\frac{5}{11}}_T+ \frac{\delta \pi \Gamma^2}{4} T^2 \|w_{HPZ}\|^{\frac{10}{11}}_T  + \frac{1}{2} \left(1 + \frac{\delta \pi \Gamma^2}{3} T^3 \right)  \|\nabla_x V(0) \|_{L^2(\Rset^3)}^2 + \|w_{HPZ} \|^{\frac{5}{11}}_T \int_0^t {\widetilde{F}}(t-s) {\widetilde{F}}'(t-s) \Big)  \\ &&+ KQ T(1+T) \\ &\leq& C(\delta, \Gamma, w_0, Q, T) \, .
%\end{eqnarray*}

%To check (v), it suffices to multiply Eq. (\ref{UZ}) (resp. (\ref{HPZeq})) by $|\xi|^2$, integrate against $x$ and $\xi$, and use the Poisson equation for the potential. 

%Assertion (vi) is a simple consequence of (\ref{ehus}) and the fact that $w_H \geq 0$.  

%Finally, (vii) follows COMPLETAR

\hfill
$\square$

Now, by taking the $L^{1,1}$ and $L^{1,2}$ norms in the mild formulation of the UZ and HPZ equations (cf. (\ref{mild})) we find

\begin{eqnarray*}
\|w(t)\|_{L^{1,1}} + \|w(t)\|_{L^{1,2}} &\leq& \|w_0\|_{L^{1,1}} + \|w_0\|_{L^{1,2}} \\ &&+ \, C \int_0^t \|H(s)\|_{L^1(\Rset^3_\xi)} \big( \|w(s)\|_{L^{1,1}} + \|w(s)\|_{L^{1,2}} \big) \, ds \\ &\leq& \|w_0\|_{L^{1,1}} + \|w_0\|_{L^{1,2}} \\ &&+ \, C \int_0^t \big(\|n(t)\|_{L^1(\Rset^3)} + \|n(t)\|_{L^2(\Rset^3)}\big) \big( \|w(s)\|_{L^{1,1}} + \|w(s)\|_{L^{1,2}} \big) \, ds \\ &\leq& \|w_0\|_{L^{1,1}} + \|w_0\|_{L^{1,2}} \\ &&+ \, C  \int_0^t  \big(Q + Q^{\frac{1}{4}} E[w](s)^{\frac{3}{4}} \big) \big( \|w(s)\|_{L^{1,1}} + \|w(s)\|_{L^{1,2}} \big) \, ds \, ,
\end{eqnarray*}
where we used Proposition \ref{apriori} (ii) and the inequality (\ref{LT}) with $p=2$.
We then conclude that
$$
\|w(t)\|_{L^{1,1}} + \|w(t)\|_{L^{1,2}} \leq \big( \|w_0\|_{L^{1,1}} + \|w_0\|_{L^{1,2}} \big) \, e^{\lambda(t)}
$$
via Proposition \ref{prop-global} (iv) and Gronwall's inequality, where $\lambda(t)$ grows less than a polynomial of at most third degree.  
%$$
%p(t) = C \int_0^t \big(Q + Q^{\frac{1}{4}} E[w](s)^{\frac{3}{4}} \big) \, ds \leq (\alpha + \beta(t)) t \, .
%$$ 
Consequently, the solutions of the UZ and HPZ initial value problems cannot blow-up at finite time, hence we are done with the proof of Theorem \ref{Th}.

%%%%%%%%%%%%%%%%%%%%%%%%%%%%%%%%%%%

\subsection*{Acknowledgements}

The first author would like to thank Departamento de Matem\'atica Aplicada of University of Granada (Spain),
where part of this work was completed, for its kind hospitality and support. He was funded by Product. CNPq
grant (Brazil) no. 305205/2016-1. The second author was supported in part by MINECO (Spain), Project MTM2014-53406-R, FEDER resources, as well as by Junta de Andaluc\'{\i}a Project P12-FQM-954.

%%%%%%%%%%%%%%%%%%%%%%%%%%%%%%%%%%%%%%%%%%%%%%%%%%%%%%%%%%%%%%%%%%%%%%%%%%%%%%

%%%%%%%%%%%%%%%%%%%%%%%%%%%%%%%%%%%%%%%%%%%%%%%%%%%%%%%%%%%%%%%%%%%%%%%%%%

\end{document}